\RequirePackage{fix-cm}
\documentclass[smallextended]{svjour3}       

\smartqed  
\usepackage{graphicx}
\usepackage{algorithmic}
\usepackage{amsmath}
\usepackage{amssymb}
\usepackage{amsfonts}
\usepackage{hyperref}
\usepackage{commutative-diagrams}
\usepackage{algorithm}

\usepackage[square,sort,comma,numbers]{natbib}

\newcommand{\Cref}[1]{Eq. \ref{#1}}
\newcommand{\cref}[1]{Fig. \ref{#1}}

\DeclareMathOperator{\diag}{diag}  

\begin{document}

\title{Visualizing fluid flows via regularized optimal mass transport with applications to neuroscience
\thanks{This study was supported by grants from the Air Force Office of Sponsored Research (FA9550-17-1-0435, FA9550-20-1-0029), a grant from National Institutes of Health (R01-AT011419).}}


\author{Xinan Chen         \and
        Anh Phong Tran \and
        Rena Elkin \and
        Helene Benveniste \and
        Allen R. Tannenbaum
}

\institute{Xinan Chen \at
              Department of Applied Mathematics \& Statistics at Stony brook University, 100 Nicolls Rd, 11794, Stony Brook, NY, USA. \\
           \and
           Anh Phong Tran $\cdot$ Rena Elkin \at
              Department of Medical Physics at Memorial Sloan Kettering Cancer Center, 1275 York Ave, 10065, New York, NY, USA.\\
           \and
           Helene Benveniste \at
              Department of Anesthesiology at Yale School of Medicine, 333 Cedar St, 06510, New Haven, CT, USA. \\
           \and
           Allen R. Tannenbaum \at
              Departments of Computer Science and Applied Mathematics \& Statistics at Stony brook University, 100 Nicolls Rd, 11794, Stony Brook, NY, USA. \\\email{allen.tannenbaum@stonybrook.edu}
}

\date{Received: date / Accepted: date}

\maketitle

\begin{abstract}
Regularized optimal mass transport (rOMT) problem adds a diffusion term to the continuity equation in the original dynamic formulation of the optimal mass transport (OMT) problem proposed by Benamou and Brenier. We show that the rOMT model serves as a powerful tool in computational fluid dynamics (CFD) for visualizing fluid flows in the glymphatic system. In the present work, we describe how to modify the previous numerical method for efficient
implementation, resulting in a significant reduction in computational runtime. Numerical results applied to
synthetic and real-data are provided.

\end{abstract}

\section{Introduction}

Optimal mass transport (OMT) treats the problem of optimally transporting a mass distribution from one configuration to another via the minimization of a given cost function. The OMT problem was first posed by Monge in 1781 in the context of the transportation of debris \cite{monge1781memoire}. This formulation was later given a modern relaxed formulation by Kantorovich in 1942 \cite{kantorovich1942translocation}. In 2000, Benamou and Brenier reformulated OMT into a computational fluid dynamics (CFD) framework \cite{benamou2000computational}. In recent times, OMT theory has received extensive research attention with rich applications in machine learning \cite{torres2021survey}, image processing/registration \cite{fitschen_rgb,peyre,haker},  network theory \cite{Buttazzo}, and biomedical science \cite{zhang2021review}.

The model employed in this work is based on the CFD approach proposed by Benamou and Brenier \cite{benamou2000computational}. Here OMT is formulated as an energy minimization problem with a partial differential equation (continuity) constraint. The continuity equation in the original version only involves advection. Our implementation includes an additional diffusion term that is of importance in our studies of glymphatic flows. This leads to the present modified formulation, which is referred as the \textit{regularized} optimal mass transport (rOMT) problem. In addition to visualizing glymphatic flow, this type of model  appears in many contexts including the Schr\"odinger bridge and entropic regularization \cite{schroedinger,entropy}. 

Here we give the formal description of the rOMT model. Given two non-negative density/mass functions $\rho_0(x)$ and $\rho_1(x)$ defined on spatial domain $\Omega \subseteq \mathbb{R}^3$ with equal total mass $\int_{\Omega}\rho_0(x)dx = \int_{\Omega}\rho_1(x)dx$, we consider the following optimization problem:

\begin{equation}\label{eq:omtenergy}
    \inf_{\rho,v}\int_0^1\int_{\Omega}\rho(t,x)\parallel v(t,x)\parallel^2dxdt
\end{equation}

subject to

\begin{subequations}
\begin{align}
  & \frac{\partial\rho}{\partial t} + \nabla\cdot(\rho v) = \sigma\Delta\rho \label{eq:omtconst1} \\
  & \rho(0,x) = \rho_0(x), \quad \rho(1,x) = \rho_1(x) \label{eq:omtconst2}
\end{align}
\end{subequations}
where $\rho(t,x): [0,1]\times\Omega\rightarrow \mathbb{R}$ and $v(t,x): [0,1]\times\Omega\rightarrow\mathbb{R}^3$ are the time-dependent density/mass function and velocity field, respectively, and $\sigma>0$ is the constant diffusion coefficient. \Cref{eq:omtconst1} is the \textit{advection-diffusion equation} or the \textit{continuity equation} in fluid dynamics. If we set $\sigma=0$, one can recover the regular OMT problem proposed by Benamou and Brenier \cite{benamou2000computational}. By adding the non-negative diffusion term $\sigma\Delta\rho$ into the continuity equation, we include both motions, advection and diffusion, into the dynamics of the system. Problem  \Cref{eq:omtenergy}-\Cref{eq:omtconst2} solves for the optimal interpolation $\rho(t,x)$ between the initial and final density/mass distributions, $\rho_0(x)$ and $\rho_1(x)$, and for the optimal velocity field $v(t,x)$ which transports $\rho_0(x)$ into $\rho_1(x)$, during which the total kinetic energy is minimized and the dynamics follow the advection-diffusion equation. Continuing the  work of the numerical method by Koundal \textit{et al.} \cite{koundal2020,elkin2018}, we report a significant reduction in runtime by about 91\% resulting from improvements of previous code.

Some of the primary applications of the present work are concerned with fluid and solute flows in the brain, and in particular, the glymphatic system. The latter is a waste clearance network in the central nervous system that is mainly active during sleep and with certain anesthetics. Many neurodegenerative diseases, such as Alzheimer’s and Parkinson’s, are believed to be related to the impairment of the function of the glymphatic system. The glymphatic transport network has received enormous attention and efforts of a number of researchers to understand the fluid behaviors in the waste disposal process in the brain \cite{nedergaard2013garbage,iliff2012paravascular,xie2013sleep,MReview2018,BReview2021}. The rOMT formulation described in the present work is highly relevant to analyzing glymphatic data due to the inclusion of both advection and diffusion terms in the continuity equation. In addition to solving the rOMT problem, we use Lagrangian coordinates for the rOMT model, which is especially useful for visualization of the time trajectories of the transport. 




\section{Material and Methods}
\label{sec:main}
This section outlines the numerical method of solving the rOMT problem \Cref{eq:omtenergy}-\Cref{eq:omtconst2}, and a post-processing Lagrangian method for practical purposes of tracing particles and visualizing fluid flows. 

\subsection{Numerical Solution of rOMT}\label{sec:numer} 
The developed method is based on the assumption that the intensity of observed dynamic contrast enhanced MRI (DCE-MRI) data is proportional to the density/mass function in the rOMT model, and thus we treat the image intensity equivalently as the concentration of the tracer molecules \textit{in vivo}. Suppose we are given the observed initial and final images, $\rho_0^{obs}(x)$ and $\rho_1^{obs}(x)$. In consideration of the image noise, instead of implementing a fixed end-point condition, we use a free end-point of the advection-diffusion process, which is realized by adding a fitting term $\parallel\rho(1,x)-\rho_1^{obs}(x)\parallel^2$ into the cost function and removing $\rho$ from the optimized variables. This free end-point version of the rOMT problem for applications in noisy (e.g., DCE-MRI) data may be expressed as
\begin{equation}
    \inf_{v}\int_0^1\int_{\Omega}\rho(t,x)\parallel v(t,x)\parallel^2dxdt + \beta\parallel\rho(1,x)-\rho_1^{obs}(x)\parallel^2 \label{eq:omtenergy_disc}
\end{equation}

subject to

\begin{subequations}
\begin{align}
  & \frac{\partial\rho}{\partial t} + \nabla\cdot(\rho v) = \sigma\Delta\rho \label{eq:omt_pde}\\
  & \rho(0,x) = \rho_0^{obs}(x)
\end{align}
\end{subequations}
where $\beta$ is the weighing parameter balancing between minimizing the kinetic energy and matching the final image. Given successive images $\rho_0^{obs},\rho_1^{obs},\cdots,\rho_{p-1}^{obs}$ where $p>2$ and $p\in\mathbb{N}^+$, this method can be recursively run between adjacent images to guide the prolonged dynamic solution.



Next, a 3D version of the algorithm is detailed. Note that the proposed workflow also works for 2D problems with simple modifications. The spacial domain $\Omega$ is
discretized into a cell-centered uniform grid of size $n_x\times n_y\times n_z$ and the time space is divided into $m$ equal intervals. Let $k_s$ and $k_t$ be the volume of each spatial voxel and the length of each time interval, respectively. With $t_i = i\cdot k_t$ for $i = 0,\cdots,m$ denoting the $m+1$ discrete time steps, we have discrete interpolations and velocity fields, $\pmb\rho = [\pmb{\rho_1};\cdots;\pmb{\rho_m}]$ and $\pmb v = [\pmb{v_0};\cdots;\pmb{v_{m-1}}]$ where $\pmb{\rho_i}$ is the interpolated image at $t = t_i$ and $\pmb v_i$ is the velocity field transporting $\pmb{\rho_i}$ to $\pmb{\rho_{i+1}}$. Note that a bold font is used to denote discretized flattened vectors to differentiate from continuous functions. For example, $\pmb{\rho}$ is a vector of size $mn\times1$ and $\pmb{v}$ is of size $3mn\times1$ where $n = n_xn_yn_z$ is the total number of voxels.

The cost function \Cref{eq:omtenergy_disc} can be approximated with
\begin{equation}
 F(\pmb v) \approx k_sk_t\pmb\rho^TM(\pmb v \odot\pmb v) +
 \parallel\pmb{\rho_m} - \pmb{\rho_1^{obs}}\parallel^2
\end{equation}
where $M=I_m\otimes[I_n|I_n|I_n]$. Here $\otimes$ is the Kronecker product; $\odot$ is the Hadamard product; $I_i$ is the $i\times i$ identity matrix; $[\cdot|\cdot]$ means forming block matrices. 

\begin{figure}
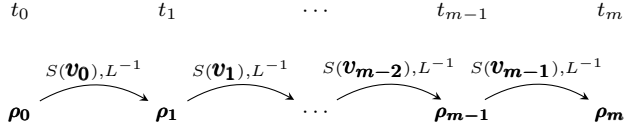

\begin{center}
\begin{codi}[l/.style={bend left}, r/.style={bend right}]
  \obj {t_0 & [2 em] t_1 & [2 em] |(pb)| \cdots & [2 em] t_{m-1} & [2 em] t_{m} \\
           |(mu0)| \pmb{\rho_0} & |(mu1)| \pmb{\rho_1} & |(pbs)| \cdots & |(mum-1)| \pmb{\rho_{m-1}} & |(mum)| \pmb{\rho_{m}} \\ };

  \mor  mu0 "\Large{S(\pmb{v_0}),L^{-1}}":l,->   mu1 "\Large{S(\pmb{v_1}),L^{-1}}":l,-> pbs "\Large{S(\pmb{v_{m-2}}),L^{-1}}":l,->  mum-1 "\Large{S(\pmb{v_{m-1}}),L^{-1}}":l,->  mum;  
\end{codi}
\end{center}
\caption{Numerical Pipeline of rOMT: From $t_i$ to $t_{i+1}$, the interpolated image $\pmb{\rho_i}$ is firstly advected via the velocity field $\pmb{v_i}$ by applying averaging matrix $S(\pmb{v_i})$ and is next diffused by applying matrix $L^{-1}$ for $i = 0,\cdots,m-1$.} \label{fig:pipeline}
\end{figure}

An operator splitting technique is employed to separate the transport process into an advective and a diffusive step. From $t_i$ to $t_{i+1}$, in the firstly advective step, a particle-in-cell method is used to re-allocate transported mass to its nearest cell centers: $\pmb{\rho_{i+1}^{adv}} = S(\pmb{v_i})\pmb{\rho_i}$ where $S(\pmb{v_i})$ is the interpolation matrix after the movement by velocity field $\pmb{v_i}$. In the secondly diffusive step, we use Euler backwards scheme to ensure stability: $\frac{1}{k_t}(\pmb{\rho_{i+1}}-\pmb{\rho_{i+1}^{adv}}) = Q\pmb{\rho_{i+1}}$ where $Q$ is the discretization matrix of the diffusive operator $\sigma\Delta$ on a cell-centered grid. Combining the two steps together, we get the discretized advection-diffusion \Cref{eq:omt_pde}:
\begin{equation}
    \pmb{\rho_{i+1}} = L^{-1}S(\pmb{v_i})\pmb{\rho_i}
\end{equation}
for $i = 0,\cdots,m-1$ where $L=I_n - k_tQ$ (\cref{fig:pipeline}). Consequently, the discrete model of the rOMT problem is given as follows:

\begin{equation}\label{eq:disc_cost}
    \min_{\pmb v}F(\pmb v) = k_sk_t\pmb\rho^TM(\pmb v \odot\pmb v) + \beta\parallel\pmb{\rho_m} - \pmb{\rho_1^{obs}}\parallel^2
\end{equation}
subject to
\begin{subequations}
\begin{align}
  \pmb{\rho_{i+1}} &= L^{-1}S(\pmb{v_i})\pmb{\rho_i}, \quad i = 0,\cdots,m-1  \label{eq:disc_pde}\\
   \pmb{\rho_0} &= \pmb{\rho_0^{obs}}.
\end{align}
\end{subequations}
One can prove that $F(\pmb v)$ is quadratic in $\pmb v$ and $S(\pmb{v_i})$ is linear in $\pmb{v_i}$. Hence, following  Steklova and Haber \cite{eldad}, the Gauss-Newton method is used to optimize for the numerical solution where the gradient $g(\pmb{v})=\nabla_{\pmb{v}}F(\pmb{v})=\frac{\partial F}{\partial\pmb{v}}$ and the Hessian matrix $H(\pmb{v})=\frac{\partial^2F}{\partial\pmb{v}^2}$ are computed to solve the linear system $Hx=-g$ for $x$ in each iteration.

Next, we elaborate on the analytical derivation of $g(\pmb{v})$ and $H(\pmb{v})$. Noticing that in $F(\pmb v)$, $\pmb{\rho}$ and $\pmb{\rho_m}$ are dependent on $\pmb{v}$ following the advection-diffusion constraint, we have 
\begin{subequations}
\begin{align}
    g = &k_sk_t\nabla_{\pmb{v}}(\pmb\rho^TM(\pmb v \odot\pmb v))+\beta\nabla_{\pmb{v}}(\parallel\pmb{\rho_m} - \pmb{\rho_1^{obs}}\parallel^2)\\
    = &k_sk_t\big(2(M\diag(\pmb v))^T\pmb\rho +  \big(\nabla_{\pmb{v}}\pmb\rho\big)^TM(\pmb v \odot\pmb v)\big) + 2\beta(\nabla_{\pmb{v}}\pmb{\rho_m})^T(\pmb{\rho_m} - \pmb{\rho_1^{obs}})
\end{align}
\end{subequations}
and
\begin{subequations}
\begin{align}
    H = \frac{\partial g}{\partial\pmb{v}} \approx&2k_sk_t\pmb\rho^T\nabla_{\pmb{v}}\big(M\diag(\pmb v)\big) + 2\beta(\nabla_{\pmb{v}}\pmb{\rho_m})^T\nabla_{\pmb{v}}\pmb{\rho_m}\\
    =&2k_sk_t\diag(\pmb\rho^TM) + 2\beta(\nabla_{\pmb{v}}\pmb{\rho_m})^T\nabla_{\pmb{v}}\pmb{\rho_m},
\end{align}
\end{subequations}
where $\diag(\cdot)$ is the function creating a diagonal matrix from the components of the given vector, and $\nabla_{\pmb{v}}$ is the operator of taking gradient with respect to $\pmb{v}$.

Considering the expressions of $g$ and $H$, the difficulty lies in the computation of $\nabla_{\pmb{v}}\pmb{\rho_m}$ and $\nabla_{\pmb{v}}\pmb{\rho}$. Let $J \triangleq \nabla_{\pmb{v}}\pmb{\rho} = (J_{\pmb{v_j}}^k)_{k,j}$ where
$J_{\pmb{v_j}}^k = \frac{\partial\pmb{\rho_k}}{\partial\pmb{v_j}}, \quad k = 1,\cdots,m, \quad j = 0,\cdots,m-1$. From the constraint \Cref{eq:disc_pde} we have
\begin{equation}\label{eq:chain}
    \pmb{\rho_k} = L^{-1}S(\pmb{v_{k-1}})L^{-1}S(\pmb{v_{k-2}})\cdots L^{-1}S(\pmb{v_{0}})\pmb{\rho_0}, \quad k = 1,\cdots,m,
\end{equation}
indicating that $\pmb{\rho_k}$ is only dependent on $\pmb{v_{0}},\cdots,\pmb{v_{k-1}}$ but independent of $\pmb{v_{k}},\cdots,\pmb{v_{m}}$, so that
for $\forall j\geqslant k, J_{\pmb{v_j}}^k=0$ holds. Therefore, $J$ is an upper-triangular block matrix of the form
\begin{equation}
    J =\begin{pmatrix}
        J_{\pmb{v_0}}^1 &  & & \\
        J_{\pmb{v_0}}^2 & J_{\pmb{v_1}}^2 & & \\
        \vdots & \vdots & \ddots &\\
        J_{\pmb{v_0}}^m & J_{\pmb{v_1}}^m & \cdots & J_{\pmb{v_{m-1}}}^m\\
        \end{pmatrix}
    \triangleq
    \begin{pmatrix}
        J_1 \\
        J_2 \\
        \vdots\\
        J_m \\
        \end{pmatrix},
\end{equation}
where $J_k = \nabla_{\pmb{v}}\pmb{\rho_k}$ is the row block of $J$ for $k = 1,\cdots,m$. If $j<k$, 
\begin{equation}\label{eq:chain2}
    J_{\pmb{v_j}}^k = L^{-1}S(\pmb{v_{k-1}})\cdots L^{-1}S(\pmb{v_{j+1}})L^{-1}B(\pmb{\rho_j}),
\end{equation}
where $B(\pmb{\rho_j})=\frac{\partial}{\partial\pmb{v_j}}(S(\pmb{v_j})\pmb{\rho_j}),$ which by the particle-in-cell method is linear to $\pmb{v_j}$ and dependent on density $\pmb{\rho_j}$. Notice that the second term of the Hessian matrix $H$ given above, involves computing the multiplication of two matrices of sizes $3mn\times n$ and $n\times 3mn,$ which is usually avoided in numerical implementation. Instead, we use a function handle that computes $Hx$ in place of the coefficient matrix $H$ so that the second term in $Hx$ can be derived from twice the multiplication of a matrix and a vector. To sum up, we are going to compute
\begin{subequations}
\begin{align}
    g = &k_sk_t\big(2(M\diag(\pmb v))^T\pmb\rho + J^TM(\pmb v \odot\pmb v)\big) + 2\beta J_m^T(\pmb{\rho_m} - \pmb{\rho_1^{obs}}),\\
    Hx = &2k_sk_t\diag(\pmb\rho^TM)x + 2\beta J_m^T J_m x.
\end{align}
\end{subequations}

We can create two functions $getJmx$ and $getJmTy$ to compute $J_mx$ and $J_m^Ty$ for any vector $x$ and $y$, respectively. Within these two functions, $J_mx$ and $J_m^Ty$ can be computed iteratively in observation of the recursive format of \Cref{eq:chain2}. Given that $J$ is upper-triangular, $J^Ty$ can be computed by recursively calling the function referred to as $getJmTy$. We can use nested function $getJmTy(getJmx(\cdot))$ to get the second term of $Hx$. However, this algorithm spends the vast majority of time (more than 90\%) solving the linear system $Hx=-g$ with the MATLAB built-in function $pcg$. We therefore modify the previous algorithm to reduce the time spent on this computational bottleneck. 

The major contributions of the present work are as follows: 
\begin{enumerate}
    \item We pre-compute all $S(\pmb{v_j})$ and $B(\pmb{\rho_k})$ for $k = 1,\cdots,m, j = 0,\cdots,m-1$ and make them inputs when calling functions $getJmx$ and $getJmTy$ to eliminate unnecessary redundant computations of advection-related matrices. 
  \item We combine nested function $getJmTy(getJmx(\cdot))$  into one function in light of the poor performance of transferring function handles in nested functions.
  \item We add the option of running the rOMT model with multiple input images $\rho_0^{obs},\rho_1^{obs},\cdots,\rho_{p-1}^{obs}$ where $p>2$ in parallel to further reduce runtime.
\end{enumerate}
Consequently, we found a significant improvement in the efficiency of solving the linear system. See Algorithm \ref{alg:algo} for the detailed process.

\begin{algorithm}
\caption{Gauss-Newton Method}
\label{alg:algo}
\begin{algorithmic}
\STATE{Load in $\pmb{\rho_0^{obs}}, \pmb{\rho_1^{obs}}$}
\STATE{$\pmb v =$ initial guess}
\FOR{$i = 1,2,\cdots, MaxIter$}
\STATE{Compute interpolations $\pmb{\rho}=$ AdvDiff$(\pmb{\rho_0^{obs}},\pmb{v})$}
\STATE{Compute $S(\pmb{v_j})$ and $B(\pmb{\rho_k})$ for $j = 0,\cdots,m-1, k = 1,\cdots,m$}
\STATE{Compute gradient $g$ and Hessian matrix function handle $Hx$}
\STATE{Solve linear system $Hx=-g$ for $x$}
\STATE{Do line search to find length $l$}
\IF{line search fails}
    \RETURN $\pmb v$
\ENDIF
\STATE{Update $\pmb{v}=\pmb{v}+lx$}
\ENDFOR
\RETURN $\pmb v$
\end{algorithmic}
\end{algorithm}

\subsection{Lagrangian coordinates}
Instead of observing the system under the usual Eulerian coordinates, one can get a Lagrangian representation of the above framework in the standard way. This is of course very useful for tracking the trajectories of particles and for investigating the characteristic patterns of fluid dynamics. This Lagrangian method has been used as a visualization method in \cite{elkin2018, koundal2020}.

Briefly, the method begins with defining the \textit{augmented velocity field} $\tilde{v} = v - \sigma\nabla\log\rho$ and putting it into the advection-diffusion equation to get a zero on the right-hand side
\begin{equation}
    \frac{\partial\rho}{\partial t} + \nabla\cdot(\rho\tilde{v}) = \sigma\Delta\rho - \sigma\nabla\cdot(\nabla\log\rho)= 0
\end{equation}
which gives a conservation form of the continuity equation \Cref{eq:omtconst1}. We apply Lagrangian coordinates $L(t,x)$ such that
\begin{subequations}
\begin{align}
  & \frac{\partial L}{\partial t} = \tilde{v}(t,L(t,x)) \label{eq:lag_coor1}\\
  & L(0,x) = x, \label{eq:lag_coor2}
\end{align}
\end{subequations}
to track the \textit{pathlines} (i.e. trajectories) of particles with the starting coordinates at $t = 0$ (\Cref{eq:lag_coor2}) and the time-varying augmented velocity field $\tilde{v}$. Along each binary pathline, the speed $s = ||v||_2$ may be calculated at each discrete time step, forming a pathline endowed with speed information which we call a \textit{speed-line}, where $||\cdot||_2$ represents the Euclidean norm. The speed-lines indicate the relative speed of the flow over time. 

This representation provides a neat way of computing certain dimensionless constants that are very popular in CFD, in particular, the P\'eclet ($Pe$) number. It has been used by several groups \cite{Mestre,Holter9894} in neuroscience to study the motion of cerebrospinal fluid (CSF) within the brain. Of special importance is the determination of regions where advection dominates or where diffusion dominates. In our model, we define  $Pe$ number as follows:
\begin{equation}
   Pe = \frac{||v||_2}{\sigma||\nabla\log\rho||_2}.
\end{equation}
This measures the ratio of advection and diffusion. Similar to speed-lines, we can compute and endow $Pe$ along the binary pathlines to form the \textit{P\'eclet-lines}. 

In order to visualize in 3D rendering, we interpolate the speed-lines and Péclet-lines into the original grid size by taking the averages of the endowed speed and $Pe$ values within the same nearest voxel to derive the smoothed \textit{speed map} and \textit{$Pe$ map}, respectively. Additionally, the directional information of the fluid flow is thereby captured by connecting the start and end points of pathlines, obtaining vectors which we refer as \textit{velocity flux vectors}. Even though they may lose intermediate path information compared to pathlines, these vectors can provide a clearer visualization of the dmovement.

The code for the Lagrangian method is available in \href{https://github.com/xinan-nancy-chen/rOMT_spdup}{https://github.com/xinan-nancy-chen/rOMT\_spdup}.

\section{Results}
\label{sec:results}
This section comes in three parts. In the first two parts, we test our methodology on self-created geometric dataset and DCE-MRI rat brain dataset, respectively. Last but not least, we compare the upgraded rOMT algorithm with the previous one on the two forementioned datasets and report a significant saving in runtime.

\begin{figure*}[t!]
\begin{center}
\setlength{\tabcolsep}{1pt}

\begin{tabular}{cccccc}
$\rho_0$ & $\rho_1$ & $\rho_2$ & $\rho_3$ & $\rho_4$ & \\
\includegraphics[width=0.18\textwidth]{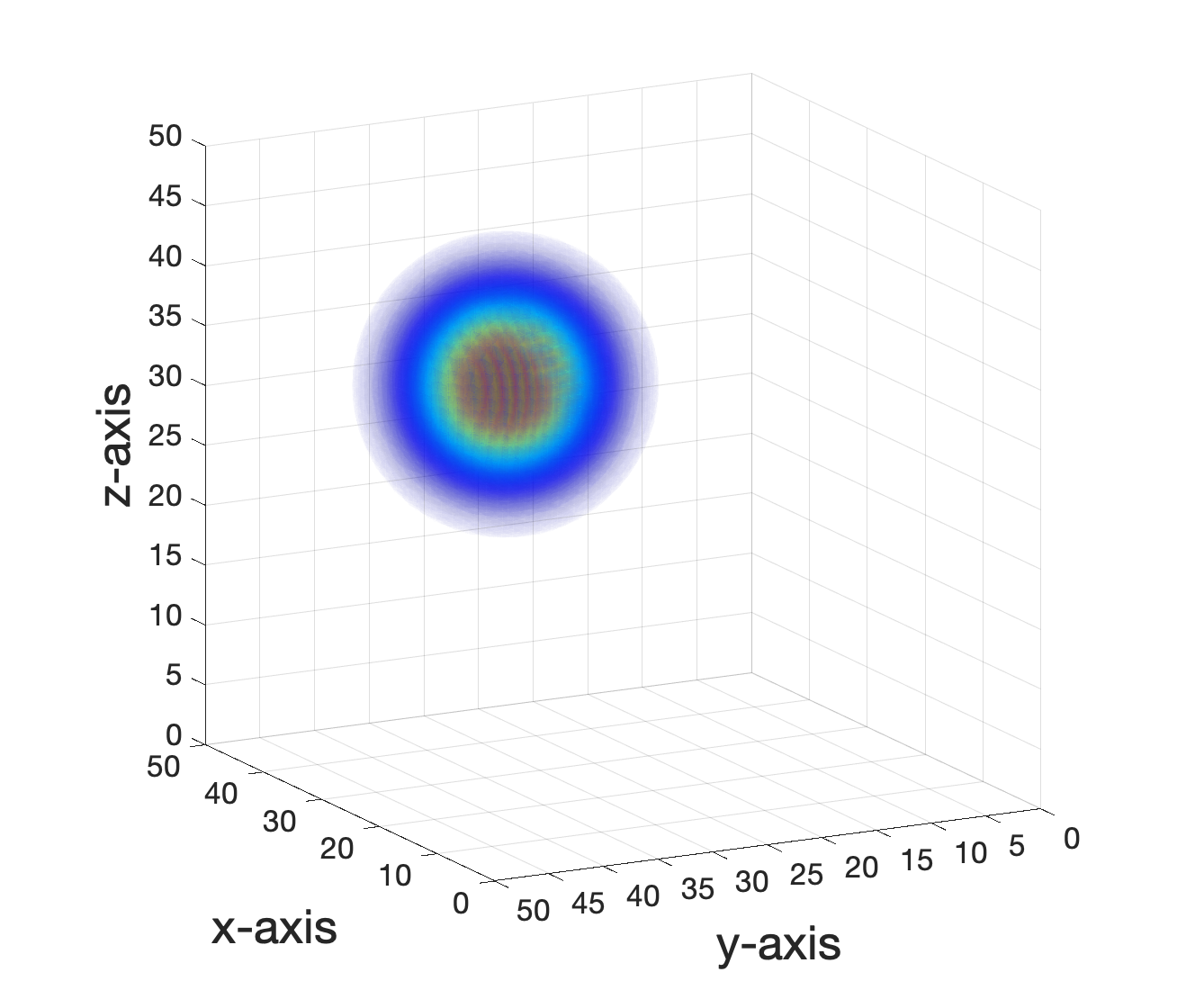}&
\includegraphics[width=0.18\textwidth]{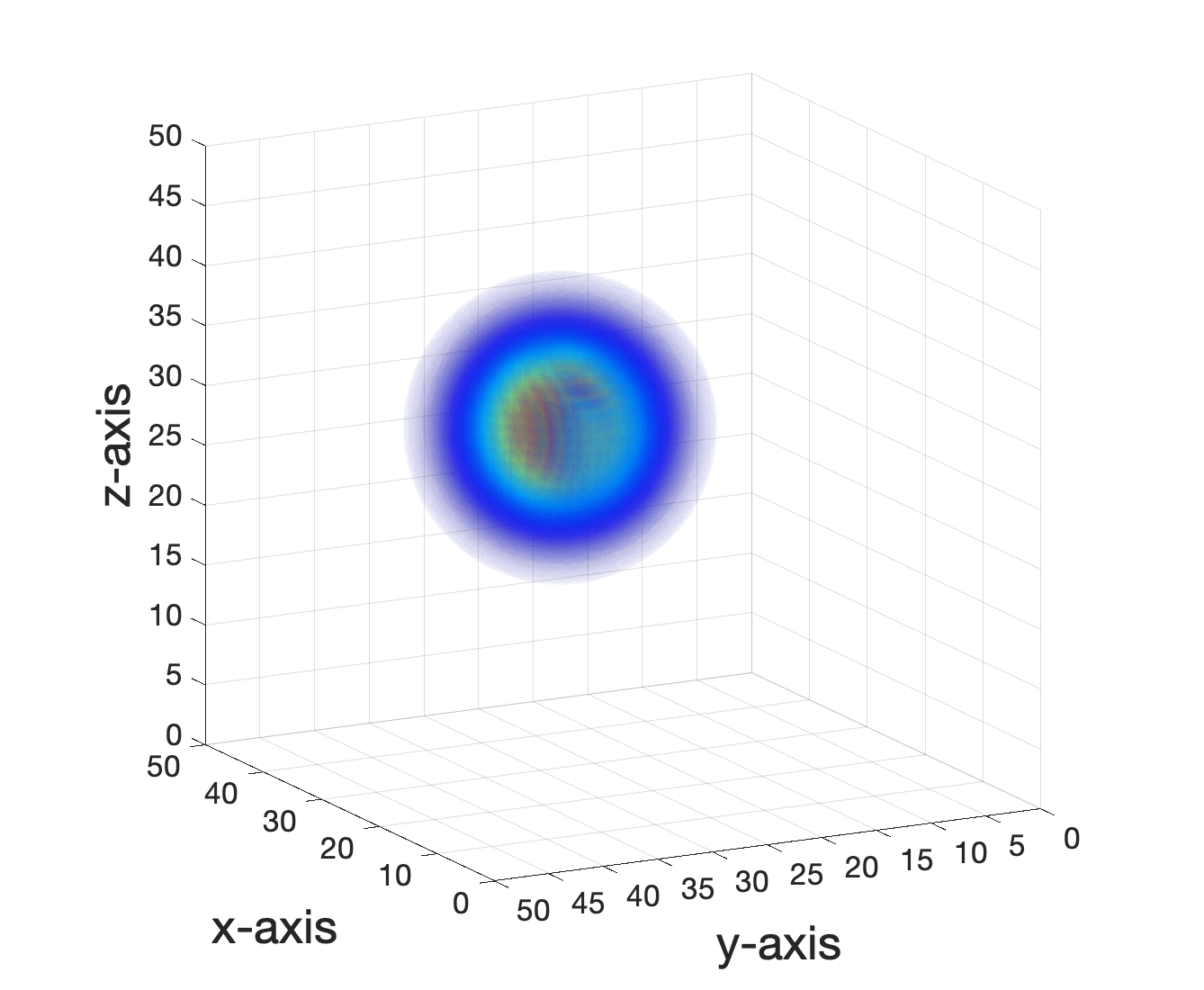}&
\includegraphics[width=0.18\textwidth]{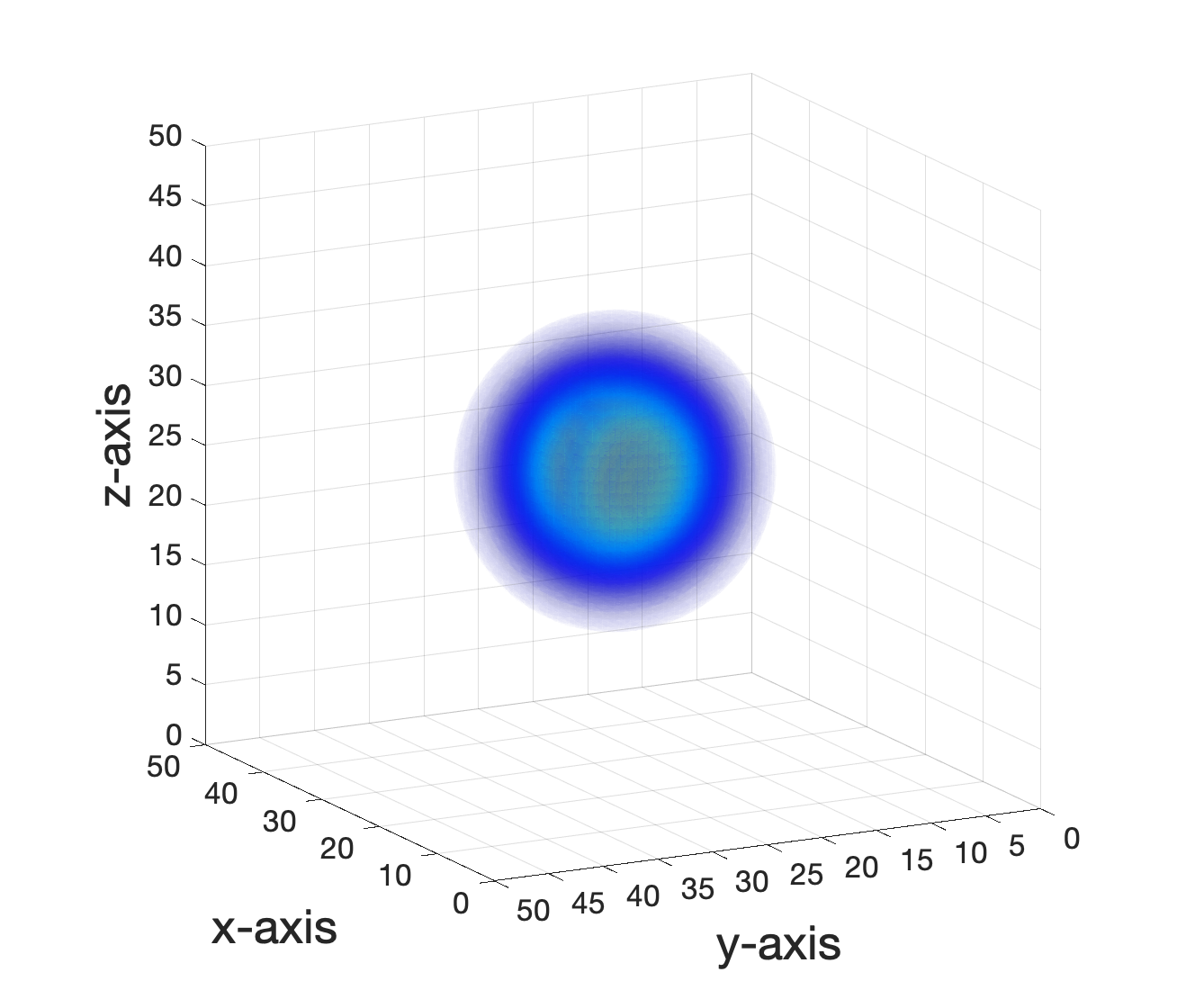}&
\includegraphics[width=0.18\textwidth]{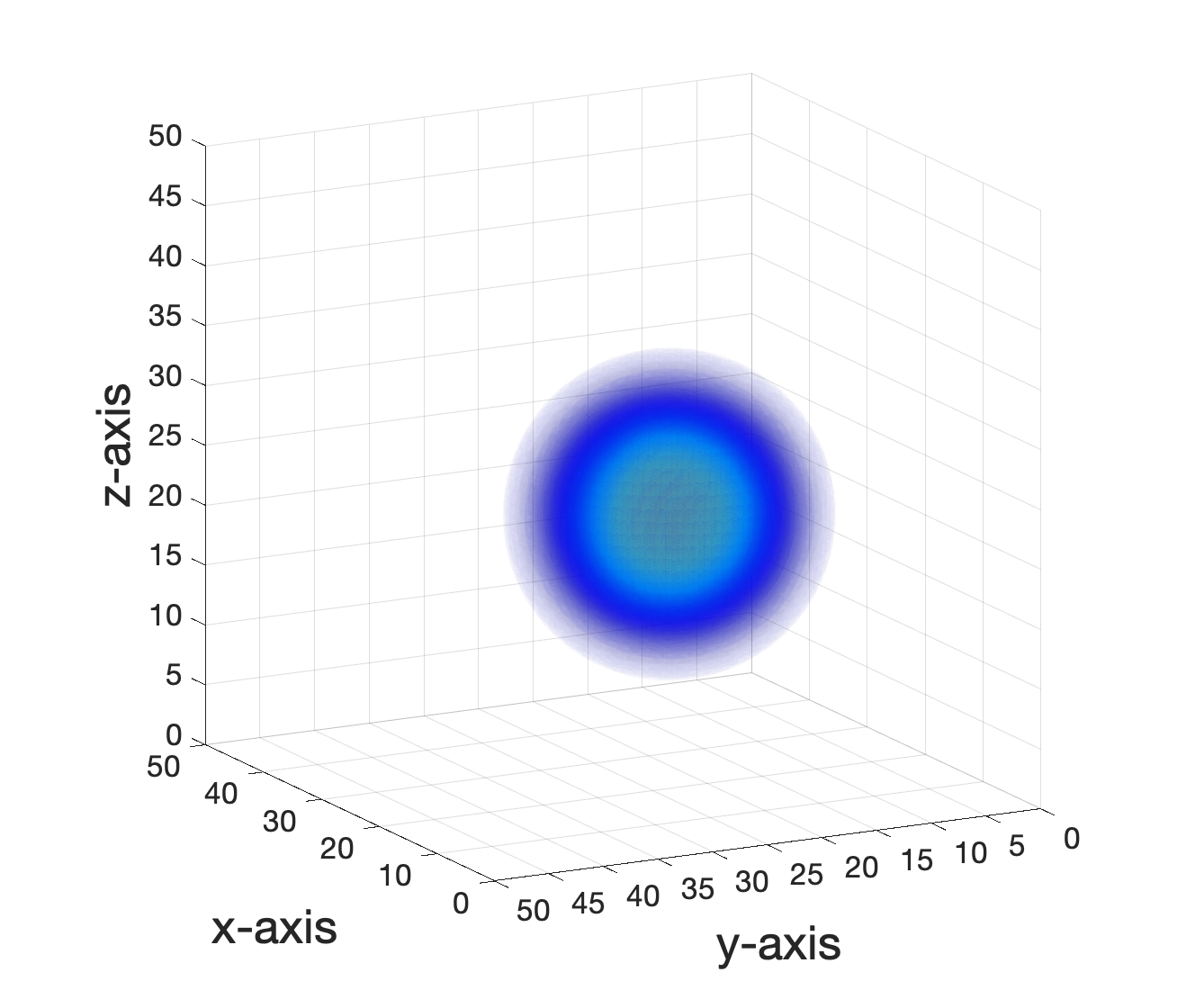}&
\includegraphics[width=0.18\textwidth]{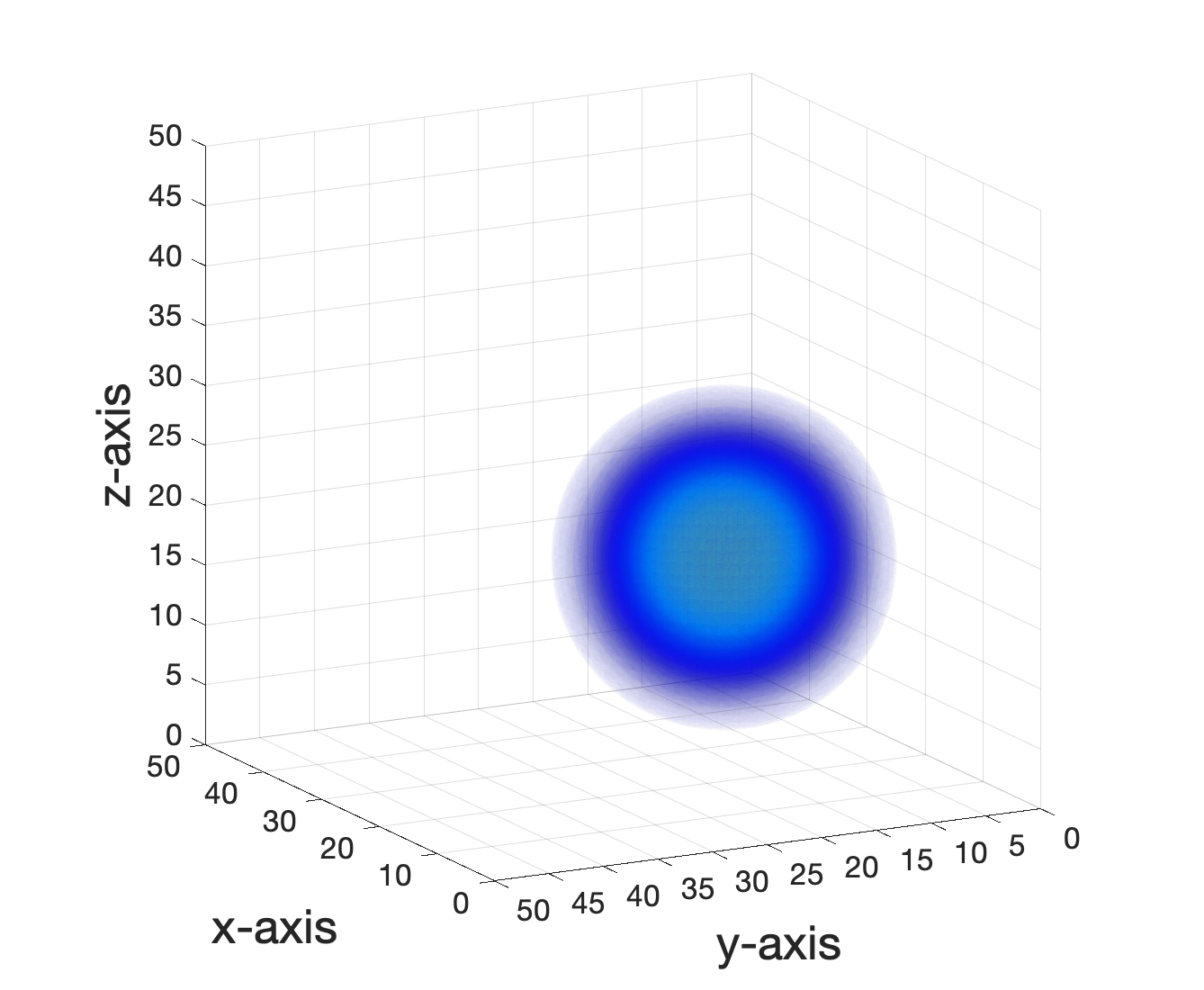}&
\includegraphics[width=0.026\textwidth]{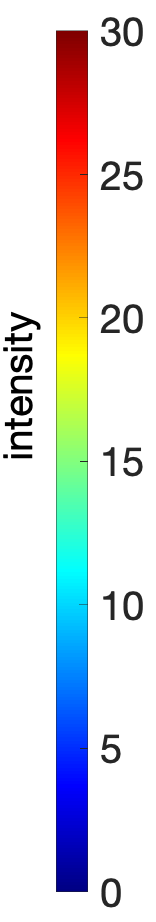}\\
\end{tabular}

\begin{tabular}{ccc}
Pathlines & Speed-lines & P\'eclet-lines\\

\includegraphics[width=0.32\textwidth]{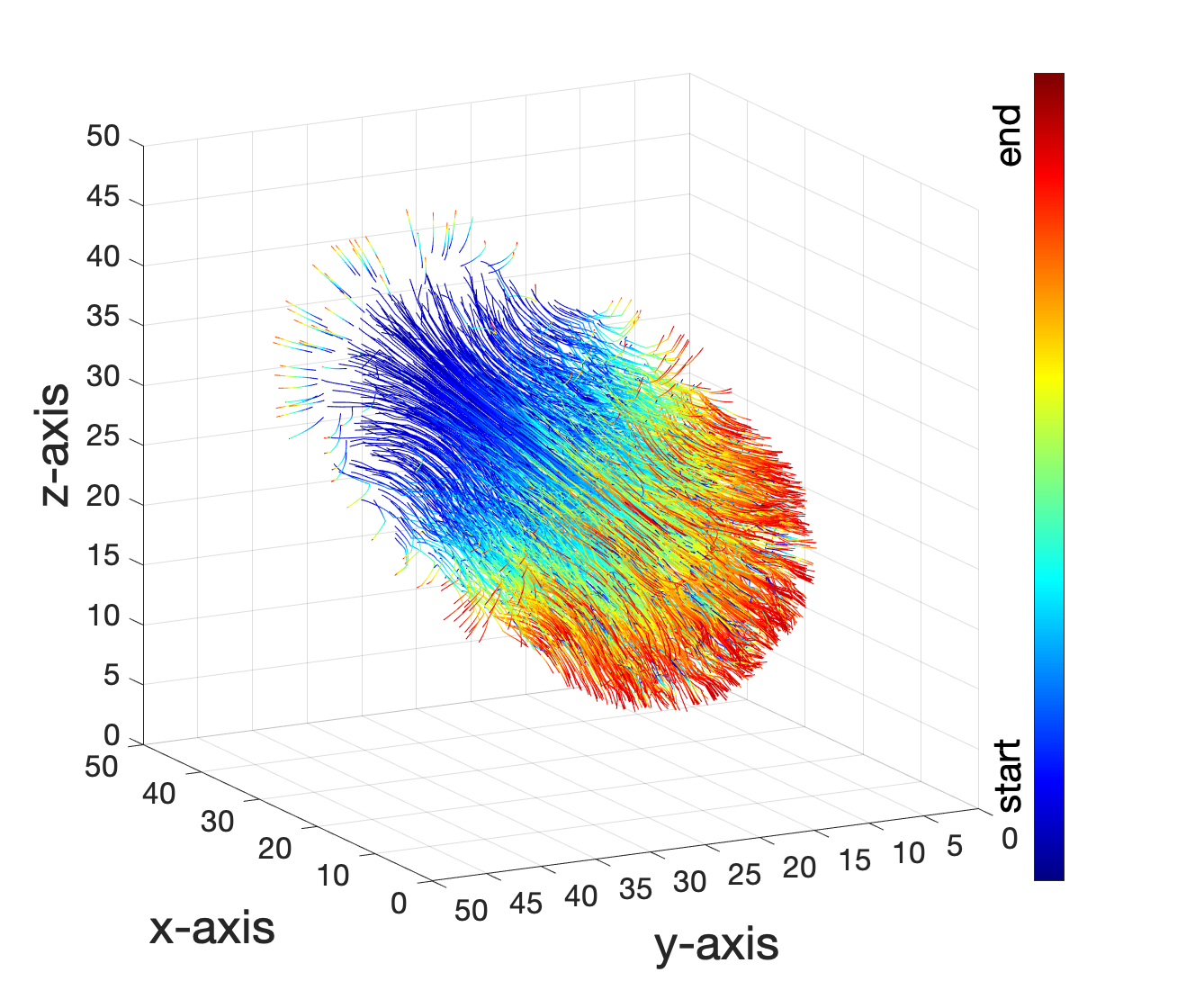}&
\includegraphics[width=0.32\textwidth]{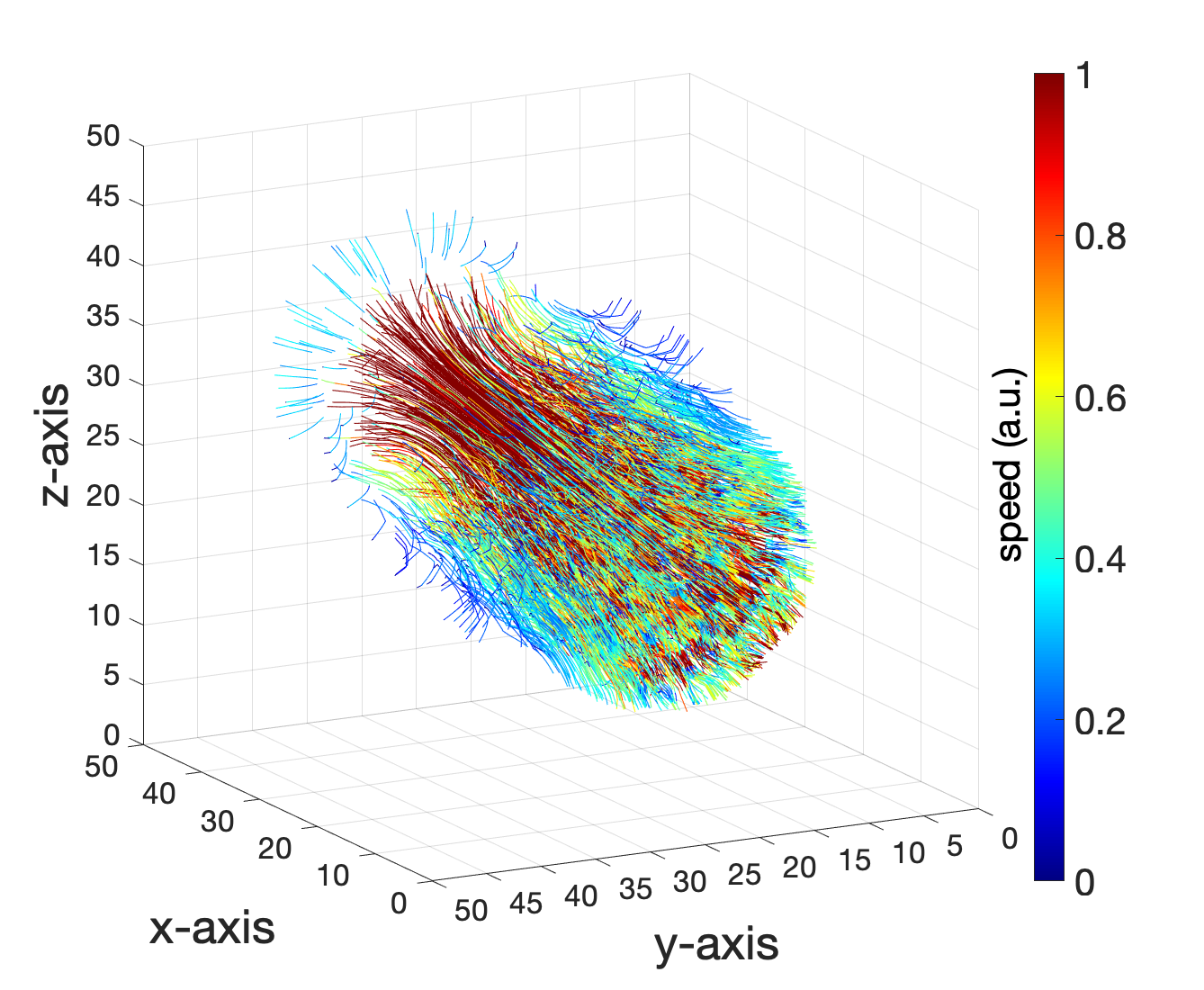}&
\includegraphics[width=0.32\textwidth]{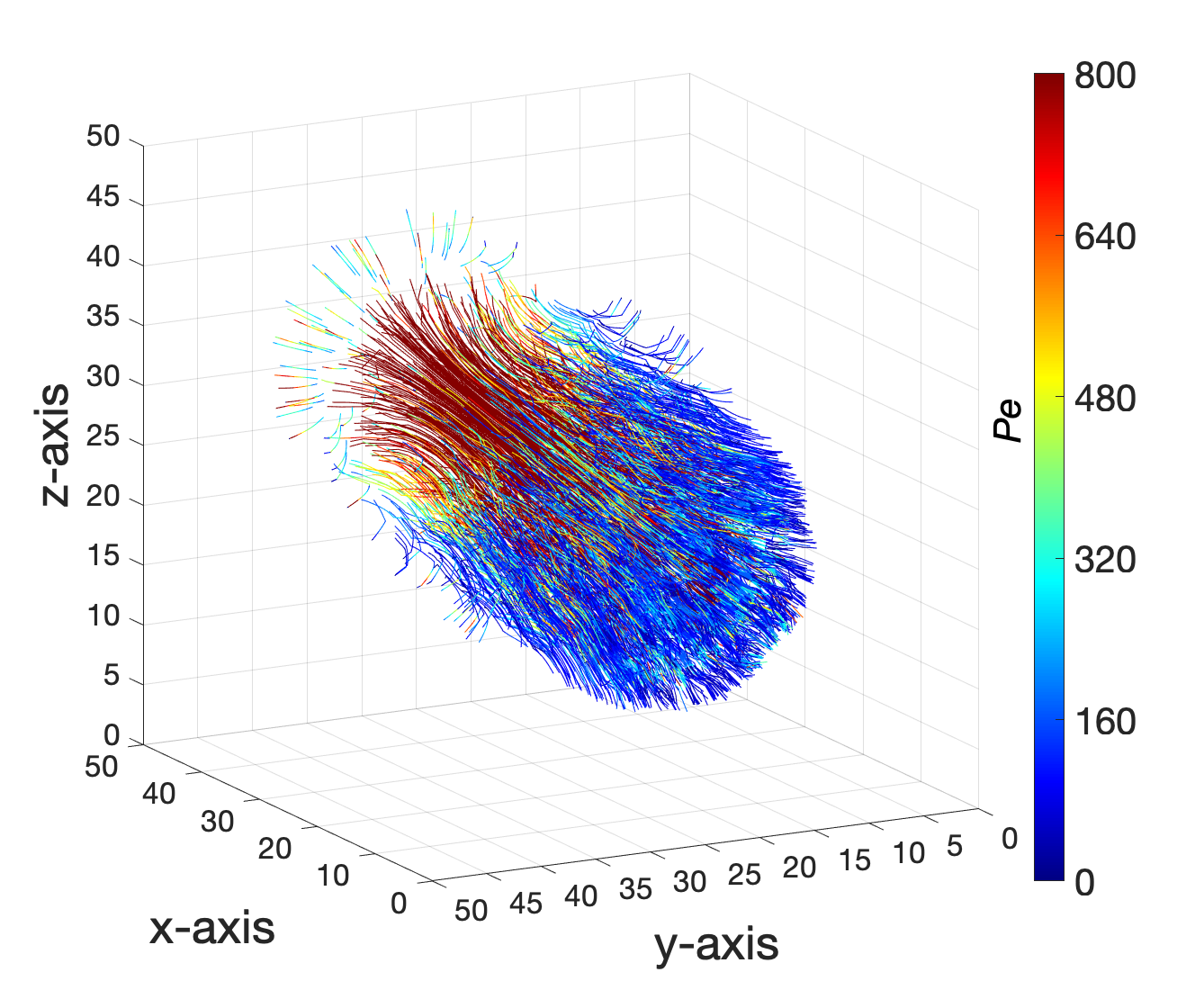}\\

\end{tabular}

\begin{tabular}{ccc}
Velocity Flux Vectors & Speed Map & $Pe$ Map\\
\includegraphics[width=0.32\textwidth]{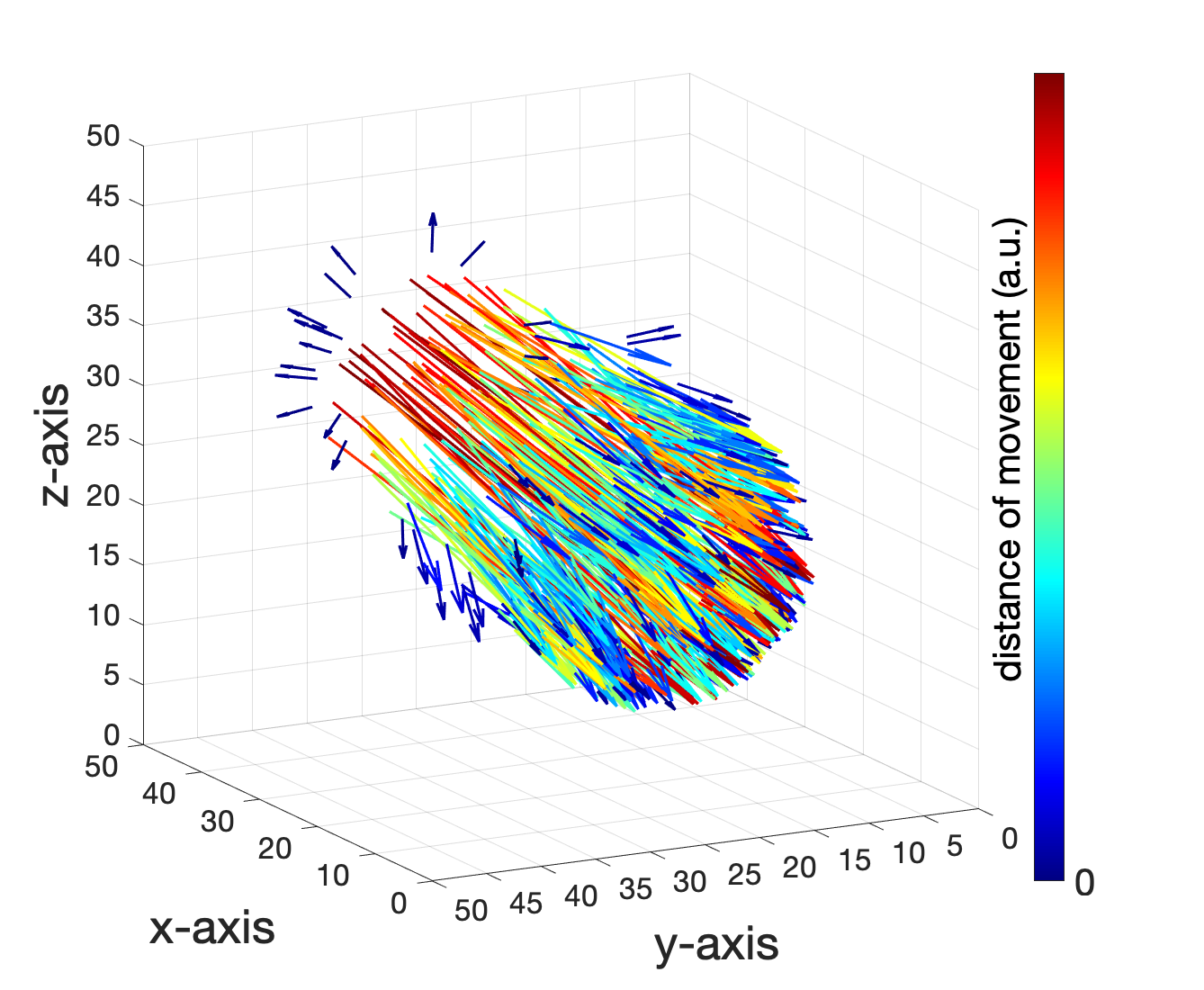}&
\includegraphics[width=0.32\textwidth]{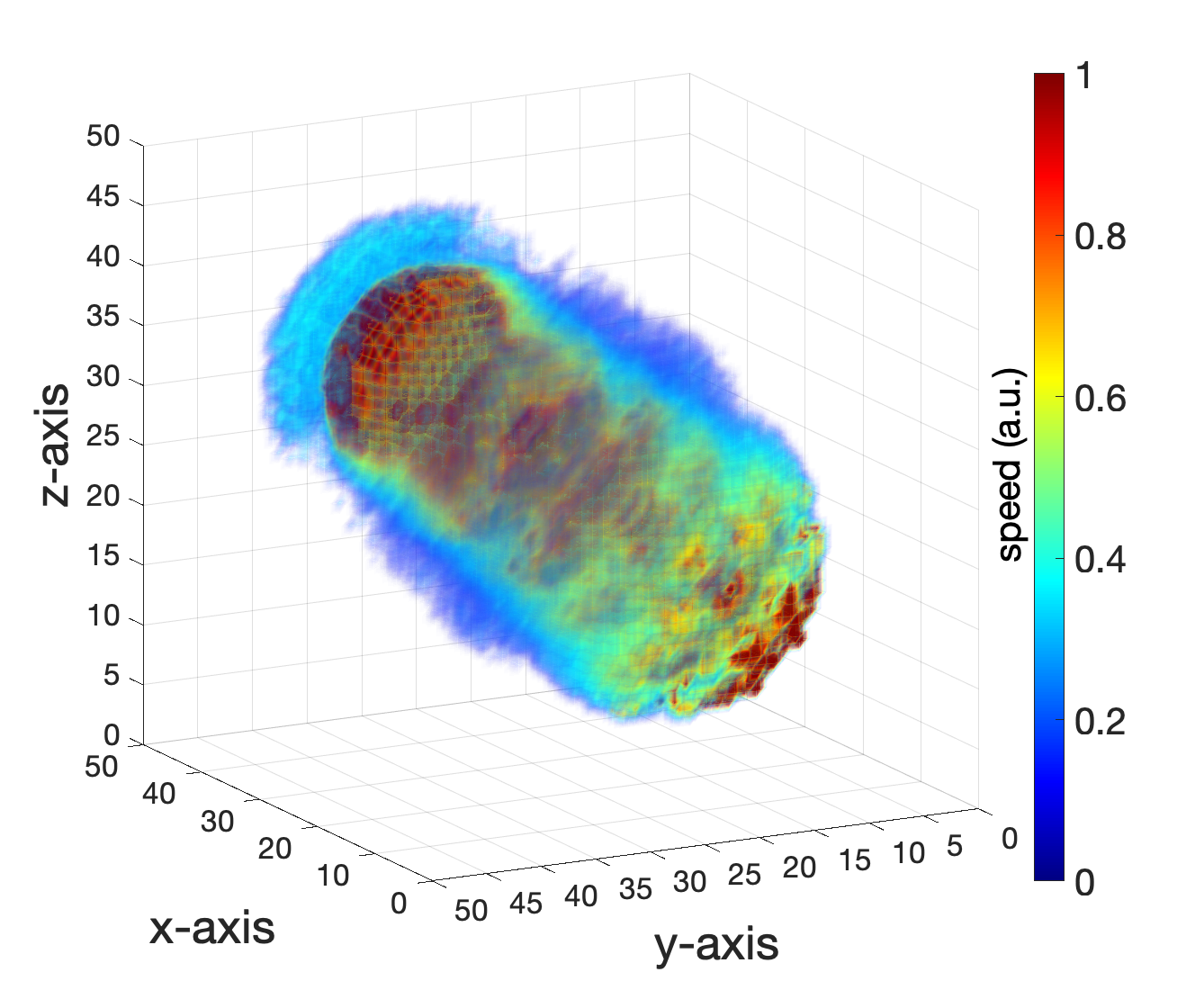}&
\includegraphics[width=0.32\textwidth]{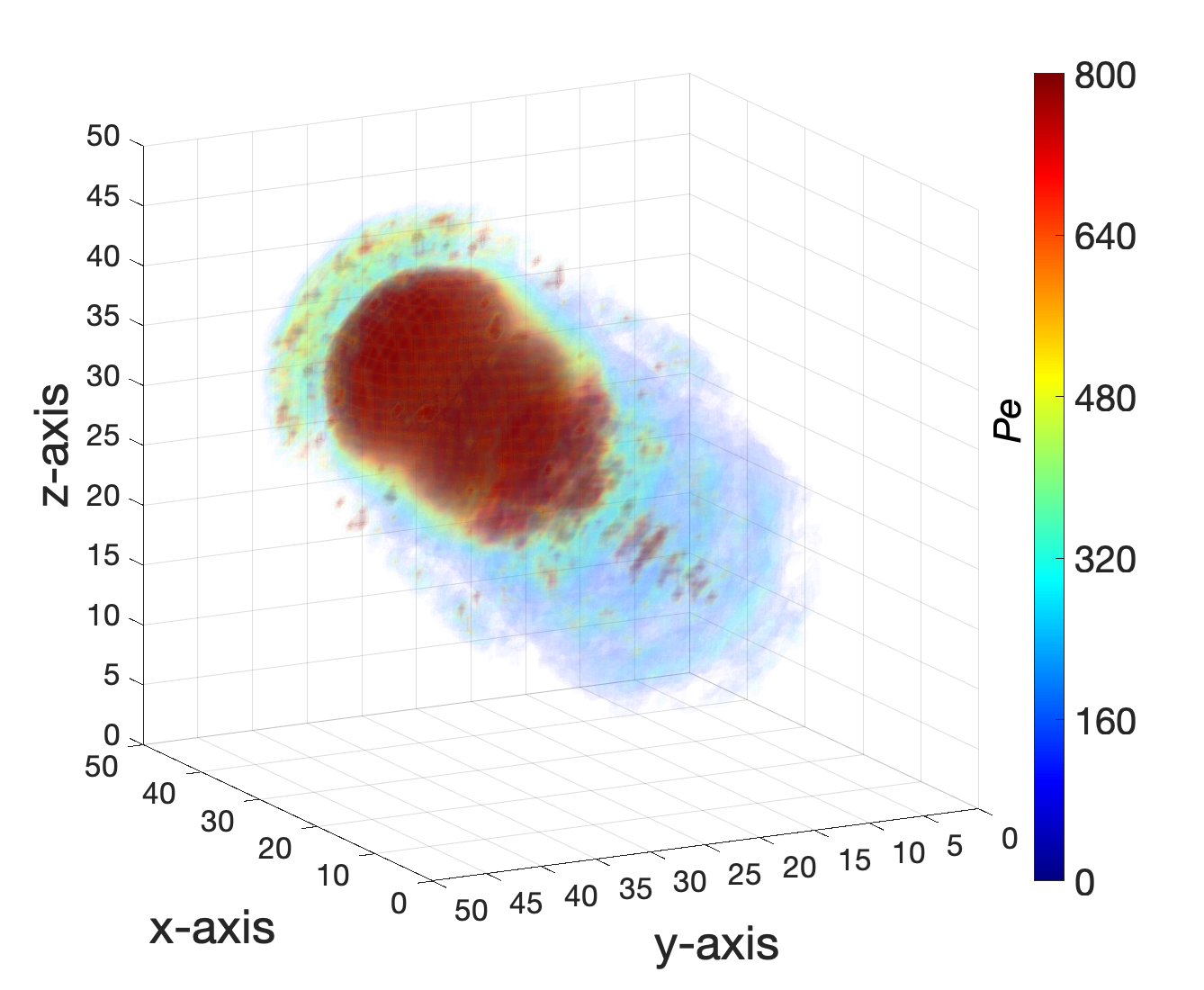}\\
\end{tabular}

\end{center}
\caption{3D Geometric Gaussian Spheres Dataset: First row: A series of synthetic Gaussian spheres as inputs of the rOMT model. The sphere is advectively moving forward while simultaneously diffusing locally. Second and third rows: Illustrative outputs from the Lagrangian  rOMT methodology to visualize the fluid dynamics. Pathlines, color-coded with start and end time, show the trajectories of particle movement. The speed-lines show the relative speed at the corresonding location along pathlines. The P\'eclet-lines indicate the local transport motion of advection-dominated (higher value) or diffusion-dominated (lower value) along pathlines. The speed map and the $Pe$ map shown in 3D rendering are smoothed interpolations of speed-lines and P\'eclet-lines on the numerical grid, respectively. Velocity flux vectors are vectors obtained by connecting the initial and terminal points of the pathlines, which are color-coded with the length of vectors, illustrating the overall direction of the movement. These outputs show that the higher speed distributes mainly in the core of the Gaussian sphere and that the transport is first dominated by advection but quickly diffusion prevails.}
\label{fig:geo_st}
\end{figure*}

\subsection{Gaussian Spheres}\label{sec:gauss}
Five 3D Gaussian spheres of image size $50\times50\times50$ are created, $\rho_0,\cdots,\rho_4$, as successive input images fed into the rOMT algorithm and its Lagrangian post-processing (\cref{fig:geo_st}, first row). The initial mass distribution is a \textit{3D} dense Gaussian sphere, and it moves forward (advection) with mass gradually diffusing into the surrounding region over time. The unparalleled runtime for this dataset was about 26 minutes on a 2.6 GHz Intel Core i7-9750H, 16 GB RAM, running macOS Mojave (version 10.14.6) with MATLAB 2019b. Please refer to Table \ref{tab:param} for parameters used in the experiment. The resulting interpolated images from the rOMT algorithm are made into a video on \href{https://github.com/xinan-nancy-chen/rOMT_spdup/blob/main/gauss2_interp.gif}{Github}. Other returned outputs are illustrated in \cref{fig:geo_st}. The binary pathlines, which are color-coded with the numerical start and end time, show the trajectories of particles. The resulting velocity flux vectors point in the direction of the movement and the color code shows relatively how far a particle is transported during the whole process. The speed-lines and the interpolated speed map indicate that the core of the Gaussian spheres are of higher speed compared to the outer regions. The P\'eclet-lines and $Pe$ map show that in the early stage of the transport, the motion is mainly advective in nature. However, in the later time, diffusion takes over. This change of dominated motion is within expectation in that dissolvable substances always have the tendency to eventually be equally mixed with the solution as a result of diffusion, regardless of the existence of an imposed velocity field.

\begin{table}[htbp]
\begin{centering}
    \begin{tabular}{ |p{1.9cm}|p{4.3cm}|p{2cm}|p{2cm}| }
    \cline{1-4}
    \textbf{Parameter} & \textbf{Definition} & \textbf{Value for Geometric Data} & \textbf{Value for Brain Data}      \\ \cline{1-4}
    \hline\hline
    $n_1$ & grid size in $x$ axis & 50 & 100 \\ \cline{1-4}
    $n_2$ & grid size in $y$ axis & 50 & 106 \\ \cline{1-4}
    $n_3$ & grid size in $z$ axis & 50 & 100 \\ \cline{1-4}
    $p$ & number of input images & 5 & 12 \\ \cline{1-4}
    $m$ & number of time intervals & \multicolumn{2}{c|}{10} \\ \cline{1-4}
    $k_t$ & length of each time interval & \multicolumn{2}{c|}{0.4} \\ \cline{1-4}
    $k_s$ & length of spatial grid & \multicolumn{2}{c|}{1}\\ \cline{1-4}
    $\sigma$ & diffusion coefficient & \multicolumn{2}{c|}{0.002} \\ \cline{1-4}
    $\beta$ &  weighting parameter in cost functional & \multicolumn{2}{c|}{5000} \\ \cline{1-4}
    \end{tabular}
\caption{Parameters used in rOMT algorithm}
\label{tab:param}
\end{centering}
\end{table}

\subsection{DCE-MRI Rat Brain}\label{sec:mri}
To further test our method for practical uses, we ran our algorithm on a DCE-MRI dataset consisting of 55 rat brains. During the MRI acquisition, all the rats were anesthetized and an amount of tracer, gadoteric acid, was injected into the CSF from the neck, moving towards the brain. The DCE-MRI data were collected every 5 minutes and were further processed to derive the \% signal change from the baseline. Data for each rat contained a 110-minute time period from 23 images of size $100\times106\times100$, and we put every other image (in total 12 images) within a masked region into our Lagrangian rOMT method to reduce the computational burden. To avoid constantly introducing new data noise into the model, we utilize the final interpolated image from the previous loop as the initial image of the next loop. The computation of the rOMT model was performed consecutively using MATLAB 2018a on the \href{https://it.stonybrook.edu/help/kb/understanding-seawulf}{Seawulf} CPU cluster using 12 threads of a Xeon Gold 6148 CPU, which took about 4 hours for each rat. The Lagrangian post-processing method took between 2 to 3 minutes for each case. Please refer to Table \ref{tab:param} for parameters used in this experiment.

\begin{figure*}[t!]
\begin{center}
\setlength{\tabcolsep}{1pt}

\begin{tabular}{c}
\includegraphics[width=0.95\textwidth]{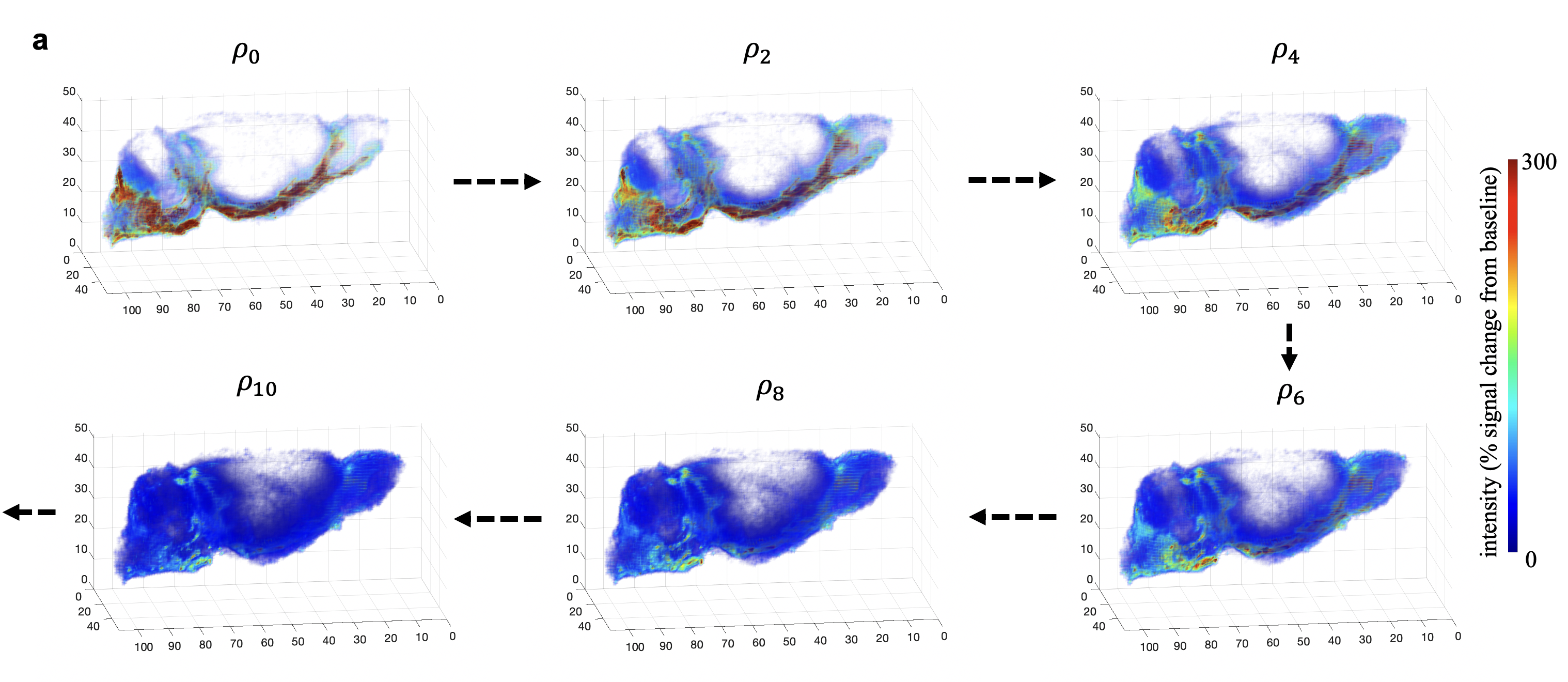}\\
\end{tabular}

\begin{tabular}{c}
\includegraphics[width=.95\textwidth]{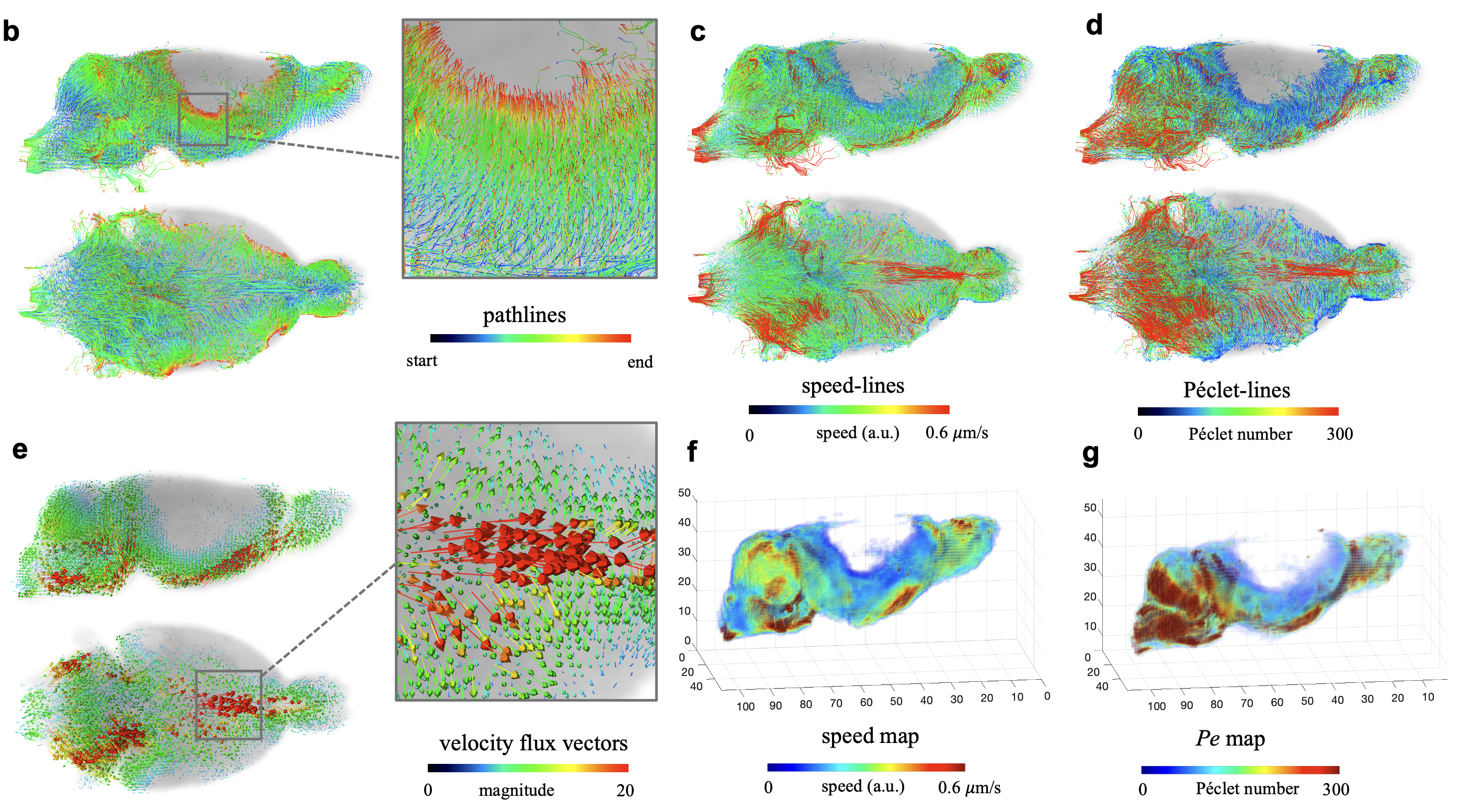}
\end{tabular}

\end{center}
\caption{3D DCE-MRI Rat Brain Dataset: (a) The input data is 12 successive images, $\rho_0,\cdots,\rho_{11}$, selectively shown in 3D rendering. (b-g) The outputs from the Lagrangian rOMT methodology. The pathlines give the trajectories of tracer over the 110-minute period. The velocity flux vectors, color-coded with the length of vectors and shown scaled by 0.25, indicate that in addition to penetrating into the brain parenchyma, there are strong flows moving within CSF towards the nose. According to the speed-lines and speed map, higher speed occurred mainly along the large vessels and quickly slowed down after entering the brain. The P\'eclet-lines and $Pe$ map identified the transport along vessels and in CSF mainly as advection-dominated by the relative higher $Pe$ numbers therein. In contrast, the movement motion was dominated by diffusion after the entry of tracer into deeper brain.}
\label{fig:mri}
\end{figure*}

In \cref{fig:mri}, we display the data and results of an example 3-month-old rat. As the pathlines and velocity flux vectors illustrate, the tracer partly entered the brain parenchyma via the CSF sink and partly was drained out towards the nose. From the speed-lines and speed map, the higher speed occurred mainly along the large vessels, which is also recognized as advection-dominated transport according to the P\'eclet-lines and $Pe$ map. When the tracer entered the brain parenchyma, the movement motion was mainly dominated by diffusion due to the relative low values there.

\subsection{Saving of Runtime}
As detailed in Section \ref{sec:numer}, we improved the current rOMT code by eliminating repeated computation in nested functions and by optimization of the algorithm. An important part of this reduction in computational time was done by pre-computing intermediate results of the advective steps that were used either throughout the calculations or for a specific step. We compared the upgraded algorithm with the previous algorithm \cite{koundal2020} by recording the runtime of rOMT code on the Gaussian sphere dataset at scaled image size (N=6, each of 4 loops) in Section \ref{sec:gauss} and the DCE-MRI dataset in Section \ref{sec:mri} (N=55, each of 11 loops). 

\begin{figure*}[t!]
\begin{center}
\setlength{\tabcolsep}{1pt}

\begin{tabular}{cc}
\includegraphics[width=0.32\textwidth]{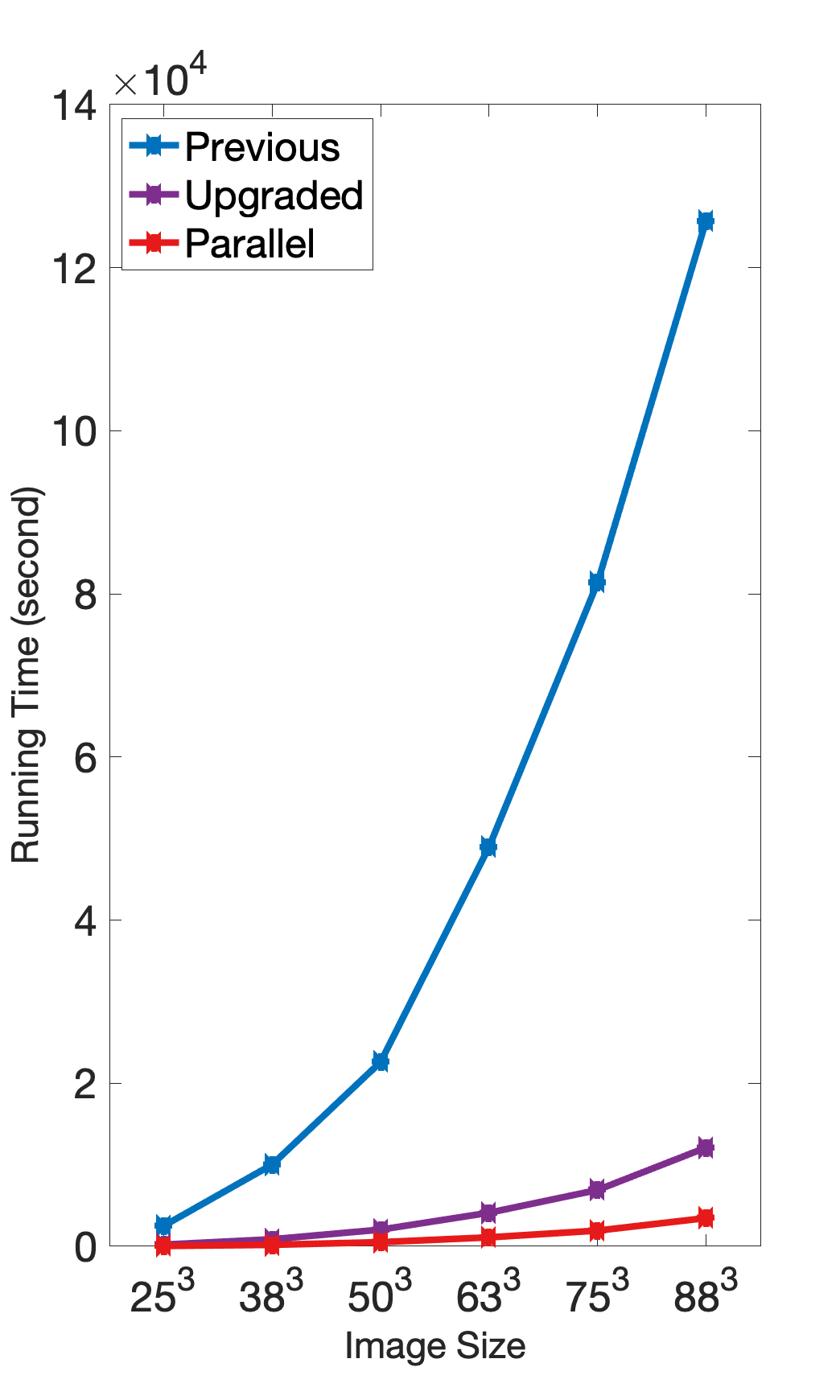}&
\includegraphics[width=0.32\textwidth]{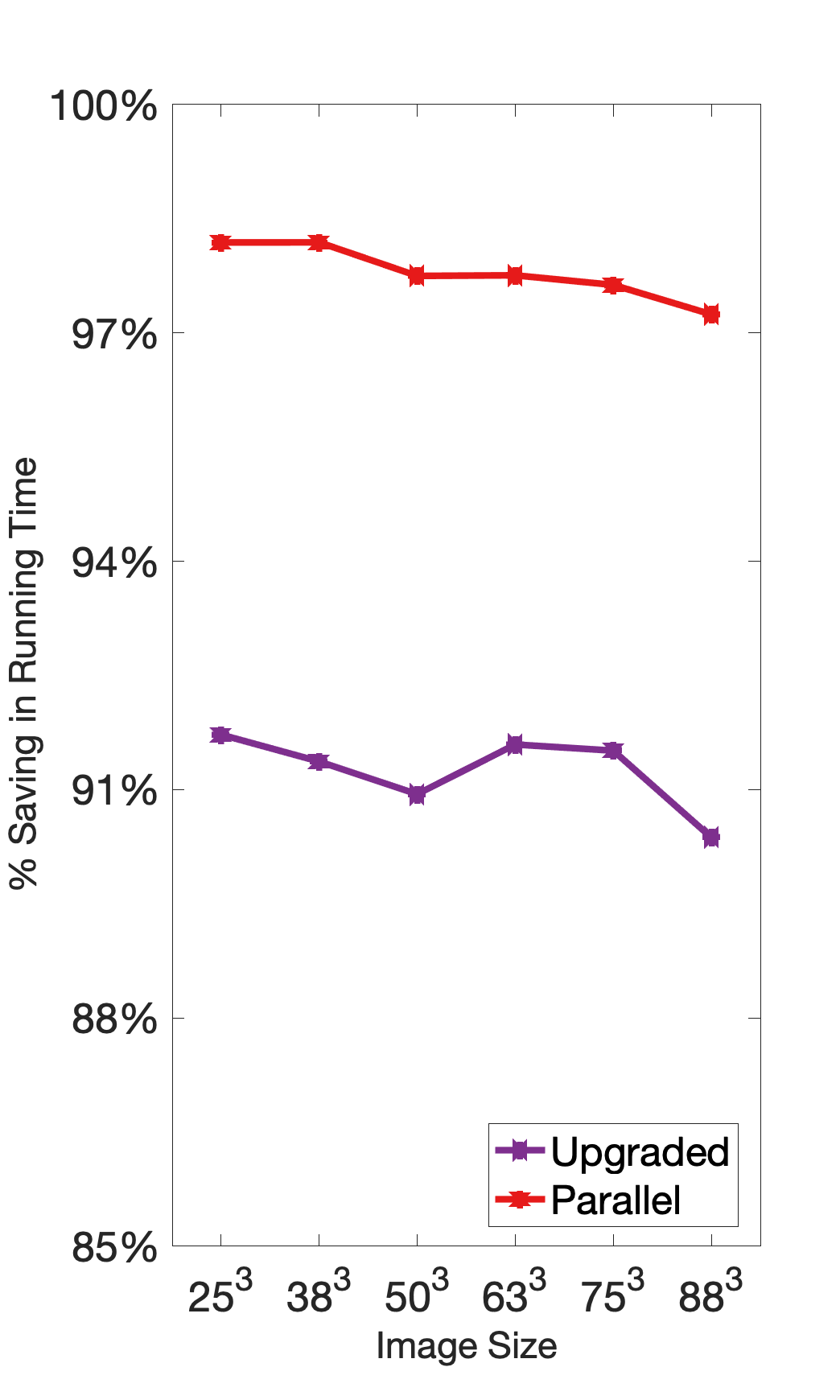}\\
\end{tabular}
\end{center}
\caption{Comparison of Runtime on the Gaussian Sphere Dataset: Left: The comparison of runtime of previous, upgraded and upgraded+parallel code at scaled image size. Right: The percent of saving in runtime compared with the previous code at scaled image size. The upgraded code realized 91.25\% $\pm$ 0.51\% reduction and can be further improved to 97.79\% $\pm$ 0.36\% if run in parallel.}
\label{fig:runningtime2}
\end{figure*}

To analyze the time complexity of the algorithm, we considered the original Gaussian sphere images of size $N^3 = 50^3$ and scaled $N$ by a factor of 0.5, 0.75, 1, 1.25, 1.5 and 1.75 to obtain spheres of sizes $25^3, 38^3, 50^3, 63^3, 75^3$ and $88^3$, respectively.  We put the six groups of data into the rOMT code using MATLAB 2019a on the \href{https://it.stonybrook.edu/help/kb/understanding-seawulf}{Seawulf} CPU cluster with 12 threads of a E5-2683v3 CPU. As illustrated in \cref{fig:runningtime2}, the runtime increases drastically as 
$N$ linearly scales up, especially for the previous code. For example, the previous code took 0.70 hours to run on the $25^3$ size input, 13.61 hours on the $63^3$ size input, and 1 day and 10.91 hours on the $88^3$ size input. The upgraded code without parallelization greatly reduced the runtime to 3 minutes, 1.14 hours, and 3.36 hours, respectively. By analyzing the six groups of data statistically, we found that in general a 91.25\% $\pm$ 0.51\% reduction and a 97.79\% $\pm$ 0.36\% reduction in runtime were realized by the upgraded code and the parallelized code, respectively, compared with the previous version.

\begin{figure*}[t!]
\begin{center}
\setlength{\tabcolsep}{1pt}

\begin{tabular}{cc}
\includegraphics[width=0.32\textwidth]{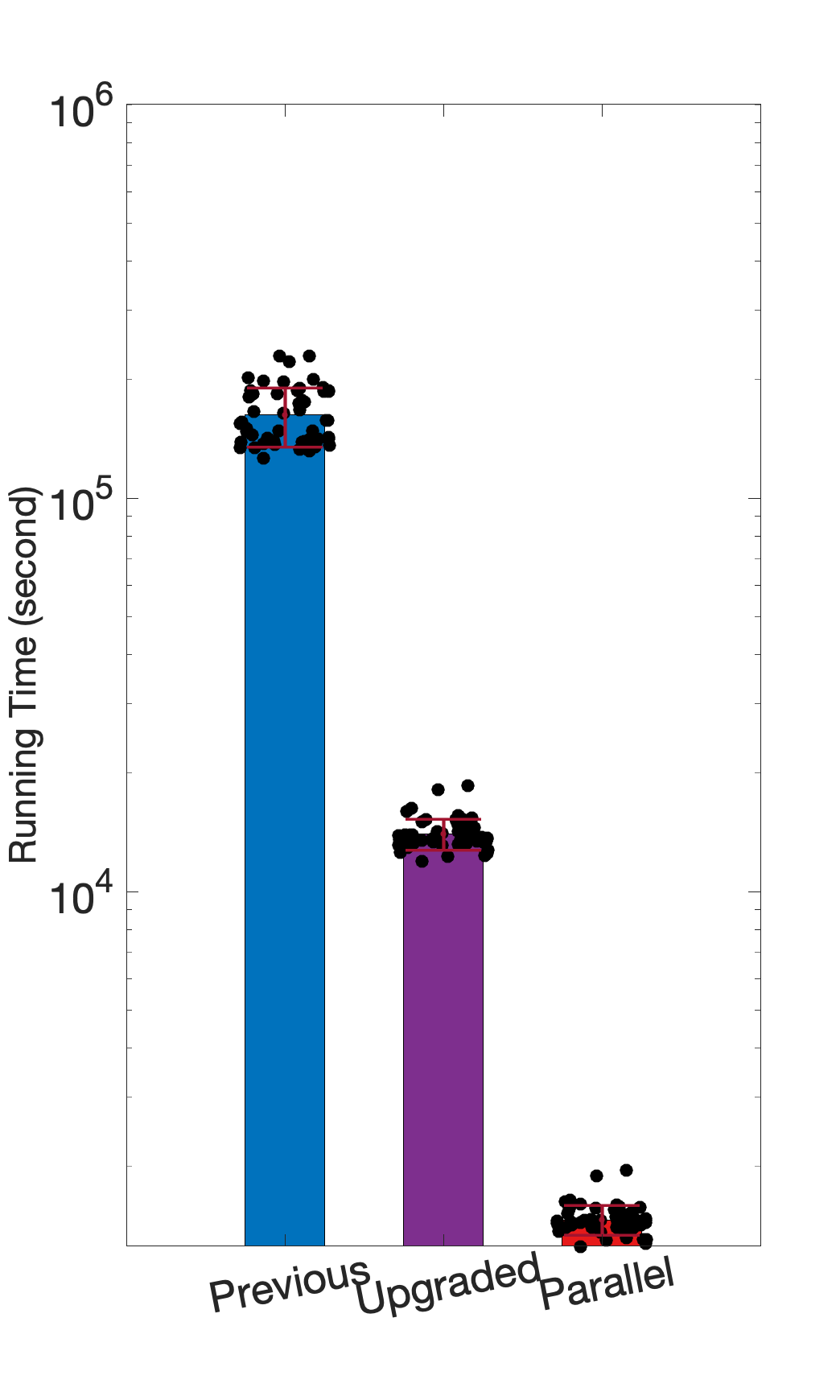}&
\includegraphics[width=0.32\textwidth]{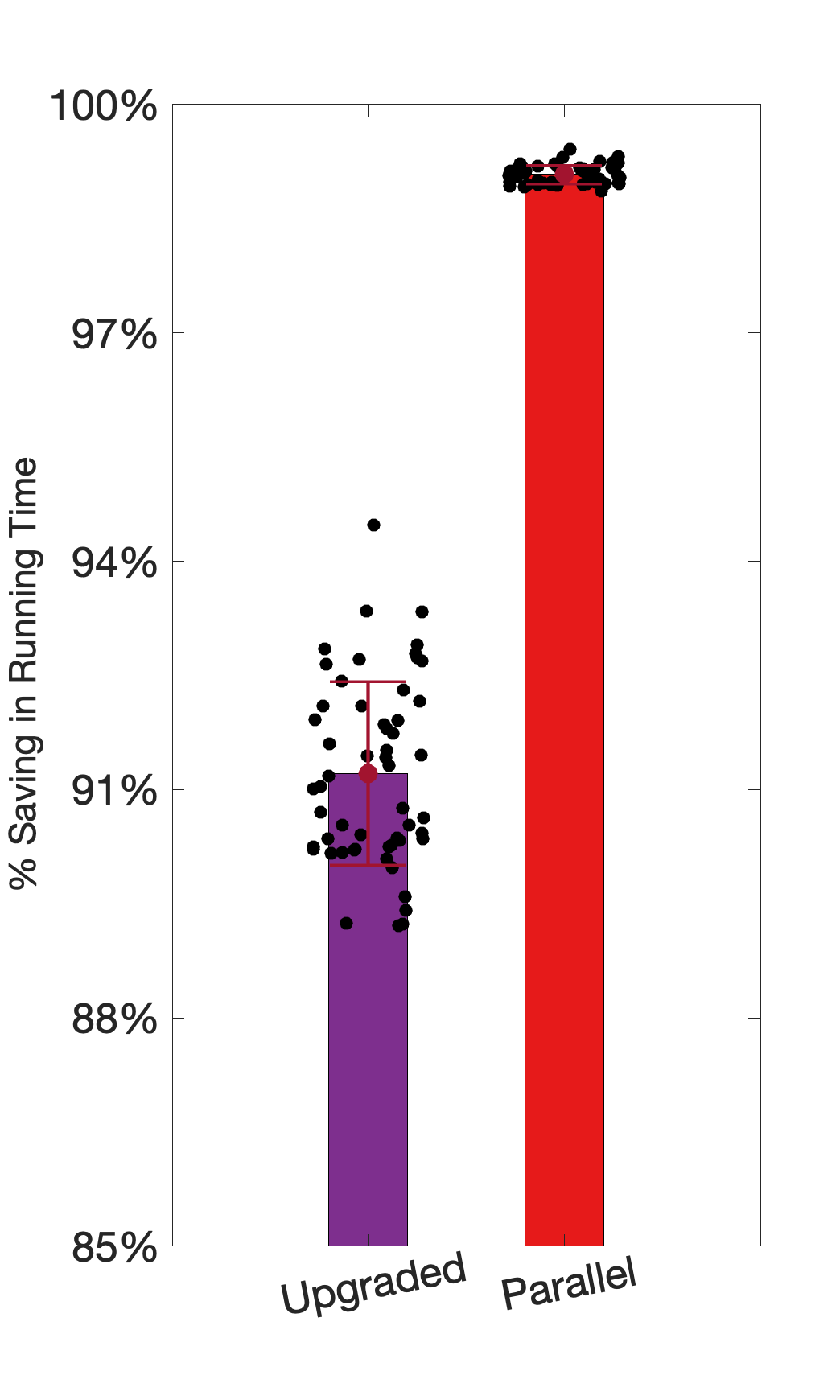}\\
\end{tabular}
\end{center}
\caption{Comparison of Runtime on the Rat Brain Dataset: Left: The comparison of runtime in logarithmic scale of previous (45.18 $\pm$ 7.64 hours), upgraded (3.89 $\pm$ 0.35 hours) and upgraded+parallel (0.41 $\pm$ 0.03 hours) code. There are in total 55 data points and each case consists of 12 input images resulting in 11 loops. Right: The percent of saving in runtime compared with the previous code. The upgraded code realized 91.21\% $\pm$ 1.21\% reduction and can be further improved to 99.08\% $\pm$ 0.12\% if run in parallel.}
\label{fig:runningtime}
\end{figure*}

For the rat brain dataset, the previous algorithm took 45.18 ($\pm$ 7.64) hours to run a case. However, it took only 3.89 ($\pm$ 0.35) hours for the upgraded algorithm to run and 0.41 ($\pm$ 0.03) hours if run in parallel, resulting in a significant reduction in runtime by 91.21\% ($\pm$ 1.21\%) and 99.08\% ($\pm$ 0.12\%), respectively (\cref{fig:runningtime}). The runtime depends on various factors, such as the size of input image, the number of input images $p$, the diffusive coefficient $\sigma$, the number of time intervals $m$, \textit{etc}. However, this significant improvement of efficiency is believed to be comprehensive as all parameters were kept fixed for the comparison.

\section{Discussion}
\label{sec:discussion}
Our rOMT methodology (both Eulerian and Lagrangian) models the dynamic fluid flows based on advection-diffusion equation and the theory of OMT. This method is largely data-driven, meaning that there is no ground truth at hand to compare with especially when it comes to real-life image data. Taking $\sigma=0$ in our model, the square root of the obtained infimum in the cost function \Cref{eq:omtenergy} gives $L^2$ {\em Wasserstein metric}, which has vast applications to many fields \cite{Villani1, Villani2}. 

It is interesting to note that rOMT is mathematically equivalent to the Schr\"odinger bridge \cite{leonard1, leonard2}, and thus the proposed algorithm may prove useful for a number of problems in which this mathematical model is relevant. We should note that
while the Schr\"odinger bridge is formally similar to OMT ($\sigma=0$), it has a substantially different motivation and interpretation. OMT was originally formulated in an engineering framework as the problem to optimally transport resources between sources and destinations. Erwin Schr\"odinger's motivation for the Schr\"odinger bridge was based on physics and the so-called ``hot gas experiment'' that led to a certain maximal likelihood problem. The aim was to link quantum theory to classical diffusion processes. In both cases (OMT and the Schr\"odinger bridge), one starts with two probability measures. In OMT, the measures are regarded as the initial and final configurations of \emph{\textbf{resources whose transportation cost has to be minimized}} among all possible couplings: this is the Kantorovich formulation \cite{Villani1, Villani2}. In comparison, the measures employed in the Schr\"odinger bridge represent initial and final probability \emph{\textbf{distributions of diffusive particles}}, and one searches for the most likely evolution from one to the other. It may be regarded as an entropy minimization problem in path space, and gives a natural data-driven model for a number of dynamical processes arising in both physics and biology, as described above in our models of fluid flow in the brain.
 
\begin{figure*}[t!]
\begin{center}
\setlength{\tabcolsep}{1pt}

\begin{tabular}{c}
\includegraphics[width=.65\textwidth]{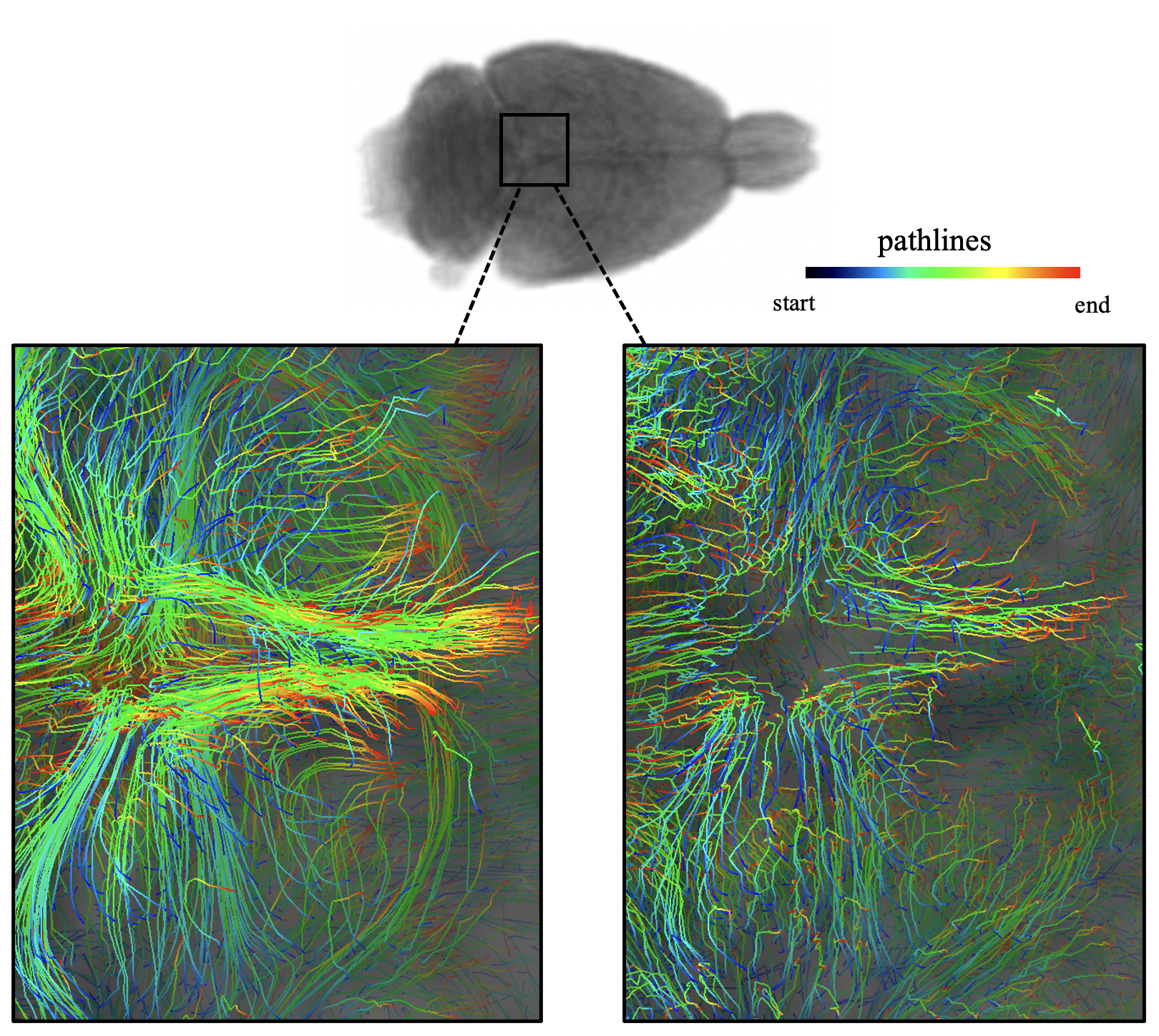}
\end{tabular}

\end{center}
\caption{Comparison of unparallelized and parallelized algorithms on pathlines: The pathlines of an example rat brain viewed from the top. Left: pathlines from continuously feeding final interpolated image into next loop. Right: pathlines from parallelized computation. The left ones are smoother and more in shape of clusters. The right ones show some fluctuations in the obtained pathlines due to constantly introducing new data noise into the system, but the overall direction of movement is aligned with the left. The left requires about 10-fold that of on the right in runtime.}
\label{fig:mri_comp}
\end{figure*}

Numerical algorithms based on finite difference scheme and optimization can be very time-consuming on large 3D images. We realized a remarkable reduction of runtime of the rOMT algorithm, cutting a two-day running down to 4 hours and even less if run in parallel. One may notice that there can be some defects in the connection of the velocity fields when applying the rOMT algorithm
independently to consecutive images
in a given series run in parallel. For example, in \cref{fig:geo_st} the speed map shows discontinued boundaries due to four independent loops. This can be ameliorated by feeding the final interpolated image from the previous loop into the next one as the initial image to give smoother velocity fields, as we did for the DCE-MRI rat brain dataset. However, by doing so parallelization is out of question because the loops are connected in time. The longer running time comes with the advantage of much smoother pathlines for meaningful visualization (\cref{fig:mri_comp}). One must be wise weighing between a quicker algorithm and smoother velocity fields. If the continuity and smoothness of velocity fields is highly emphasized, a multi-marginal model should be considered and longer time to run is also expected.



\section{Conclusions}
\label{sec:conclusions}
We introduced the rOMT methodology, both in the Eulerian and Lagrangian frameworks, as an efficient algorithm to track and visualize fluid flows, and in particular to track the trajectories of substances in glymphatic system using DCE-MRIs from a computational fluid dynamic perspective. Quantitative measurements, speed and the P\'eclet number are also provided along with the pathways to help uncover features of the fluid flows. We improved the previous code by removing redundant computation and significantly saved the runtime by 91\%, and we offer the option to further cut the runtime down by putting the algorithm in parallel.

\begin{acknowledgements}
Acknowledgements This work was supported by grants from the National Institutes of Health/
National Institute on Aging (AG053991) and the AFOSR (FA9550-20-1-0029). 
\end{acknowledgements}

%
\section*{Conflict of interest}
The authors declare that they have no conflict of interest.

\bibliographystyle{spmpsci}      
\bibliography{references}   

%
%

\end{document}